\documentclass[11pt,article]{article}
\usepackage{amsfonts,amssymb,amsmath,array} 
\usepackage[all]{xy}
\newtheorem{Lemma}{Lemma}[section]\newcommand{\bel}{\begin{Lemma}}\newcommand{\eel}{\end{Lemma}}
\newtheorem{Proposition}[Lemma]{Proposition}\newcommand{\bprop}{\begin{Proposition}}\newcommand{\eprop}{\end{Proposition}}
\newtheorem{Theorem}[Lemma]{Theorem}\newcommand{\bthe}{\begin{Theorem}}\newcommand{\ethe}{\end{Theorem}}
\newcommand{\bpr}{{\bf proof:~}}\newcommand{\epr}{$\blacksquare$\\}
\newtheorem{Remark}[Lemma]{Remark}\newcommand{\beR}{\begin{Remark}\rm}\newcommand{\eeR}{\end{Remark}}
\newtheorem{Definition}[Lemma]{Definition}\newcommand{\bed}{\begin{Definition}}\newcommand{\eed}{\end{Definition}}
\newtheorem{Example}[Lemma]{Example}\newcommand{\bex}{\begin{Example}\rm}\newcommand{\eex}{\end{Example}}
\newtheorem{Corollary}[Lemma]{Corollary}\newcommand{\bcor}{\begin{Corollary}\rm}\newcommand{\ecor}{\end{Corollary}}
\newtheorem{Fact}[Lemma]{Fact}\newcommand{\bfact}{\begin{Fact}\rm}\newcommand{\efact}{\end{Fact}}

\newcommand{\beq}{\begin{equation}}\newcommand{\eeq}{\end{equation}}
\newcommand{\bem}{\begin{displaymath}}\newcommand{\eem}{\end{displaymath}}
\newcommand{\beqa}{\begin{eqnarray}}\newcommand{\eeqa}{\end{eqnarray}}
\newcommand{\bee}{\begin{enumerate}}\newcommand{\eee}{\end{enumerate}}
\newcommand{\bei}{\begin{itemize}}\newcommand{\eei}{\end{itemize}}
\newcommand{\bet}{\begin{tabular}{cccccccc}}\newcommand{\eet}{\end{tabular}}
\newcommand{\bpm}{\begin{pmatrix}}\newcommand{\epm}{\end{pmatrix}}
\newcommand{\ber}{\begin{array}{l}}\newcommand{\eer}{\end{array}}

\newcommand{\di}{\partial}
\newcommand{\ra}{\!\!\rightarrow\!\!}

\newcommand{\tinyM}{\scriptstyle}
\newcommand{\tinyT}{\scriptsize}
\newcommand{\tinyA}{\tinyM\tinyT}

\newcommand{\cL}{{\cal{L}}}\newcommand{\cO}{{\cal{O}}}

\newcommand{\mC}{\mathbb{C}}
\newcommand{\mN}{\mathbb{N}}
\newcommand{\mP}{\mathbb{P}}\newcommand{\mS}{\mathbb{S}}\newcommand{\mXL}{\mP_x^2\times(\mP_l^2)^*}
\newcommand{\mPN}{\mP_f^{N_d}}\newcommand{\mPJ}{\mP Jet}\newcommand{\mR}{\mathbb{R}}
\newcommand{\mZ}{\mathbb{Z}}

\newcommand{\al}{\alpha}\newcommand{\be}{\beta}\newcommand{\De}{\Delta}

\newcommand{\ep}{\epsilon}
\newcommand{\Si}{\Sigma}

\newcommand{\tl}{\tilde{l}}\newcommand{\tSi}{\widetilde{\Sigma}}
\newcommand{\tv}{\tilde{v}}
\newcommand{\tpp}{\mP^2_x\underset{\De}{\tilde{\times}}\mP^2_y}

\newcommand{\hPn}{\widehat{\mP^n}}\newcommand{\hPt}{\widehat{\mP^2}}\newcommand{\hPN}{\widehat{\mP_f^{N_d}}}
\newcommand{\hPJ}{\widehat{\mPJ}}\newcommand{\hSi}{\widehat\Si}\newcommand{\htSi}{\widehat{{\widetilde\Si}}}

\newcommand{\bSi}{\bar\Sigma}\newcommand{\btSi}{\bar{\widetilde{\Sigma}}}
\newcommand{\lSi}{\boldsymbol{\overline{\Sigma}}}\newcommand{\ltSi}{\boldsymbol{\overline{\widetilde\Sigma}}}

\newcommand{\li}{~\\ $\bullet$ }
\newcommand{\into}{\stackrel{i}{\hookrightarrow}}

\newcommand{\gNN}{generalized Newton-non-degenerate }

\newcommand{\mesh}[7]
{
\put(#1,#2){\vector(1,0){#6}}  \put(#1,#2){\vector(0,1){#7}}
\setcounter{tempx}{#3}  \addtocounter{tempx}{1}  \setcounter{tempy}{#4}  \addtocounter{tempy}{1}
\multiput(#1,#2)(#5,0){\value{tempx}}{\multiput(-1.5,-0.5)(0,#5){\value{tempy}}{.}}
}

\title{On the geometry of some strata of uni-singular curves}
\author{D. Kerner
}

\begin{document}\setcounter{secnumdepth}{6} \setcounter{tocdepth}{2}\newcounter{tempx}\newcounter{tempy}
\maketitle
\begin{abstract}
We study geometric properties of linear strata of uni-singular curves. We resolve the singularities of
closures of the strata and represent the resolutions as projective bundles. This enables us to study their
geometry. In particular we calculate the Picard groups of the strata and
the intersection rings of the closures of the strata. The rational equivalence classes of
some geometric cycles on the strata are calculated.
As an application we give an example when the proper stratum is not affine.
\\\\
As an auxiliary problem we discuss the collision of two singular points, restrictions on
possible resulting singularity types and solve the collision problem in several cases. Then we present some cases of
enumeration of curves with two singular points, one of them being a node.
\end{abstract}
\tableofcontents

\section{Introduction}
\subsection{Formulation and results}
We work with (complex) algebraic curves in $\mP^2$. Identify the complete linear system $|dL|$
(the parameter space of plane curves of degree $d$) with the projective space $\mPN$.
Here \mbox{$N_d={d+2\choose{2}}-1$}, the subscript $f$ is due to the defining equation of the curves, $f(x)=0$.

The parameter space is stratified according to the {\it embedded topological singularity type} of curves.
The generic point of $\mPN$
corresponds to a smooth curve. The set of points corresponding to singular curves is called {\it the discriminant}
 ($\Sigma$). It is a (projective) hypersurface in the parameter space.

An {\it equisingular stratum} is the (quasi-projective) variety of points corresponding to the curves with the given
topological type of singularity. The generic point of the discriminant lies in the stratum of nodal curves
($\Sigma_{A_1}$). Other strata correspond to higher singularities (e.g. curves with r-nodes $\Si_{(A_1)^r}$
or $\Si_{A_k}$, $\Si_{D_k}$, $\Si_{E_k}$ etc..). For a comprehensive introduction cf. \cite{GLSbook}.
In this paper we study some of the strata and their compactifications.
To be clear we sometimes call the strata themselves: {\it the proper strata}.

The degree of curves, $d$, is assumed to be sufficiently high (for a specified singularity type).
Then the proper strata possess good geometric properties:
are non-empty, irreducible, reduced, smooth algebraic varieties of expected dimension
(for recent review cf. \cite{GLS06}).
A sufficient condition for this is:

if the curve has $r$ singularities of types $\mS_i$, with
orders of determinacy $o.d.(\mS_i)$, then the degree must be not less than $\sum o.d.(\mS_i)+r-1$
\cite[proposition I.3.9]{Dim}.
\\\\\\
The compactified strata are singular in codimension one.
In the previous work \cite{Ker06} the (partial) resolutions of the compactified strata of uni-singular
curves were constructed
 as subvarieties
of some multi-projective spaces. This enabled us to calculate their cohomology classes $[\Si_\mS]$ (in the
cohomology ring of the ambient space).

Here we study the geometry of the strata.
We restrict the consideration to curves with one singular point of {\it linear} type (the precise definition
is in section \ref{SecSingularityTypes}). The simplest examples
of linear singularities are $A_{k\le3}$, $D_{k\le6}$, $E_{k\le8}$, $X_9$, $J_{10}$, $Z_{k\le13}$ etc.
(In this paper the low codimension singularity types are denoted according to the tables in \cite[section II.15]{AGV}).
We discuss compactifications and partial resolutions of the linear strata.
\\\\\\
 The above cohomology classes are used to obtain the information
about rational equivalence on the strata. In particular we calculate
the intersection rings of the (resolved compactifications of the) strata. As an example we give explicit formulae
in the following cases (the type is specified by the normal form):
\\ $x^{p+1}_1+x^{p+1}_2$ (ordinary multiple point), $x^{p}_1+x^{p+1}_2$ (generalized cusp, e.g. $A_2,E_6..$),
$x^{p}_1+x_1x^{p}_2+x^{p+2}_2$ (e.g. $A_3,E_7..$), $x^{p+1}_1+x^2_1x^{p-1}_2+x^{p+2}_2$ (e.g. $A_2,D_5..$),
$x^{p+1}_1+x^2_1x^{p-1}_2+x^{p+3}_2$ for $p\ge3$.

We study the boundary divisors of the compactification $\Si\subset\lSi$. Their equivalence classes
were mostly obtained in the previous work, now they are used to calculate the Picard group
 of the proper strata $\Si$ or their partial compactifications. As a by-product, one can often check whether
 the proper stratum is affine (cf. \ref{SecIntroAffineStratum}).
 In particular we present examples of {\it non-affine proper strata}.

One component of the boundary always corresponds to curves with {\it two} singular points (one of them being a node).
To study this component we describe a method to enumerate curves with two singular point.
As an example we apply it in several cases: $\Si_{x^{p+1}_1+x^{p+1}_2,A_1}$, $\Si_{x^{p}_1+x^{p+1}_2,A_1}$,
$\Si_{x^{p}_1+x_1x^{p}_2+x^{p+2}_2,A_1}$, $\Si_{x^{p+1}_1+x^2_1x^{p-1}_2+x^{p+2}_2,A_1}$,
$\Si_{x^{p+1}_1+x^2_1x^{p-1}_2+x^{p+3}_2,A_1}$.  In some of these cases we are able also to calculate
the intersection ring of the (resolved compactification of the) corresponding stratum.
\subsubsection{Historic overview.} While enumerative questions and questions about irreducibility, dimension,
smoothness of the strata are very old (with results starting from the 19'th century, cf. \cite{GLS06,Kaz4}),
the questions on the intersection rings of the strata seems to be almost untouched (to the best of authors knowledge).
\li In \cite{DiazHarris86,DiazHarris88} the geometry of the Severi varieties ($\Si_{(A_1)^r},\lSi_{(A_1)^r}$)
was studied and many divisors on it were described. It was proved that the proper stratum is affine and
conjectures on the structure of $Pic(\Si_{(A_1)^r})$ and $Pic(\lSi_{(A_1)^r})$ were formulated.
For example: $Pic(\Si_{(A_1})^r)$ is torsion.
\li In \cite{MiretXamboDescamp94} the Picard groups of $\Si_{A_1}$ and $\lSi_{A_1}$ were calculated (verifying
in particular the above conjecture).
\li In \cite{MiretValls} the result was generalized to the stratum of plane
curves with ordinary multiple point $\lSi_{x^{p+1}_1+x^{p+1}_2}$. They identified the boundary divisors of
$\Si_{x^{p+1}_1+x^{p+1}_2}\subset\lSi_{x^{p+1}_1+x^{p+1}_2}$ and calculated the intersection rings $A^*(\lSi_{x^{p+1}_1+x^{p+1}_2})$ and
$A^*(\lSi_{x^{p+1}_1+x^{p+1}_2,A_1})$. In particular they verified the above conjecture for
 $Pic(\Si_{x^{p+1}_1+x^{p+1}_2})$.
\li \cite{Edidin94} has showed that $Pic(\Si_{(A_1)^r})$ is torsion for sufficiently many nodes
(for a given degree of curve $d$).
\subsubsection{On compactification and resolution.}\label{SecIntroCompactificationAndResolution}
The equisingular strata are defined as subvarieties of the parameter space $\mPN$.  Correspondingly, they
have a natural compactification: the topological closure (denoted by $\lSi$). This compactification is highly
singular (in codimension 1). Another unpleasant feature is that that the complement $\lSi\setminus\Si$
is in general not of pure dimension (so, it is not a divisor). A simple example of this is:
$\di\lSi_{A_2}=\lSi_{A_3}\cup\lSi_{A_2A_1}\cup\lSi_{D_4}$.

A natural way to resolve the singularities is to consider the {\it universal curve} (i.e. the
lifting to a bigger ambient space)
\beq\label{EqMinimalLifting}
\ltSi(x)=\overline{\left\{(x,f)|\text{The~curve}~\{f(x)=0\}~\text{has the~prescribed singularity at~the point }x\right\}}
\into\mP^2_x\times\mPN
\eeq
The projection $\ltSi(x)\to\lSi$ is surjective and generically 1:1.
In the simplest case of ordinary multiple point this "lifted" stratum is already a smooth variety. In general
one must lift further and consider "generalized universal curves", taking into account other parameter
of the singular germ: tangent lines ($l$), osculating conics etc.

We denote lifted strata by $\tSi,\ltSi$ and usually assign the parameters of lifting (as in (\ref{EqMinimalLifting})).
So, the lifted strata are subvarieties of some multi-projective spaces ($Aux\times\mPN$ with $Aux$
for the auxiliary).

For linear singularities the second lifting suffices: $\ltSi(x,l)$ is already a smooth variety. For
non-linear singularities the situation is more complicated. The (naturally) lifted strata can be
still singular and it is not clear whether they can be desingularized by lifting only.

When considering curves with several singularities, the corresponding universal curve (the minimal lifting)
is of course $\tSi(x,y,..)$: the variety of singular curves with singular points assigned.

Somewhat unexpectedly, the lifting often solves the problem of the boundary also:
\bprop\label{ThmBoundaryIsHypersurface}
For the lifted stratum of linear singularity $\ltSi(x,l)$ the boundary
$\ltSi(x,l)\setminus\tSi(x,l)$ is of pure co-dimension 1
(i.e. a true hypersurface).
\eprop
This is proved in section \ref{SecBoundaryComponents} by describing the nearest adjacent types for linear singularities.
~~\\\\
For various purposes partial compactifications are important. In particular it
is natural to consider the
compactification with the preserved topological type of the {\it chosen} singularity.
We call it the {\bf semi-compactification} and denote $\btSi,\bSi$).
So, the curve that belongs to $\bSi$ can be singular
at other points (e.g. $\bSi_{A_k}$ contains $\bSi_{A_k,A_1}$ but is disjoint to $\Si_{A_{k+1}}$).
\subsubsection{On rational equivalence and the intersection rings.}\label{SecIntroIntersectionRings}
Having constructed a lifted stratum (a smooth projective variety)
one can study various divisors/cycles on it.
\\
There are several types of divisors/cycles in $\ltSi$ (in the spirit of \cite{DiazHarris88}):
\li {\bf Boundary divisors/cycles.} Among them are:

-those corresponding to the curves with additional singular point (i.e. the strata $\ltSi_{\mS,A_1}\subset\ltSi_\mS$).

-those corresponding to degenerations of the singularity type that do not increase the degree of determinacy.

-those corresponding to degenerations of the singularity type that increase the degree of determinacy
\li {\bf Intrinsic divisors.} Various divisors/cycles of the ambient space (a multi-projective space) pulled-back
to the lifted stratum. For example, for the minimal lifting (as in (\ref{EqMinimalLifting})), let $X,F$ be the classes
dual to the corresponding hyperplanes in $\mP_x^2\times\mPN$. Then the cycle $i^*(X^jF^k)$ corresponds to the
family of plane curves (with singularity of a prescribed type) that pass through $k$ fixed (generic)
points and whose singularity lies on a fixed generic $2-i$ plane.
\li {\bf Extrinsic divisors/cycles} are defined by the properties of (embedded) curves. For example
cycles of curves with hyperflexes or multi-tangents, singular points whose smooth branches have flexes or tangent
lines are also tangent at other points etc..
\\\\\\
We are interested in (rational equivalence) classes of the divisors/cycles. It appears that classes of all
the cycles are expressible through the classes of the intrinsic divisors
(the precise statement is \ref{ThmIntroIntersRingLinearStrat}).

For boundary divisors the corresponding expression is either obtained directly
 (using the methods of \cite[Appendix A]{Ker06}) or one first should
 calculate the cohomology class of a (lifted) boundary divisor ($\ltSi_{*,A_1}\subset Aux\times\mPN$) and
 then to represent it as a product $[\ltSi_*][D]$.

To express an extrinsic divisor/cycle through the classes of intrinsic one should make additional calculations
(e.g. to impose the conditions of tangency and to intersect the lifted stratum with the corresponding hypersurfaces).
Here we consider only the cycle of curves whose smooth branches have (hyper)flexes at the singular point.
\\\\
As was noticed above, intrinsic divisors generate the intersection ring. More precisely:
\bprop\label{ThmIntroIntersRingLinearStrat}
\li For linear singularities the strata lifted to $\mXL\times\mPN$ are smooth projective bundles
(over a smooth base $\{(x,l)|~l(x)=0\}$). They are rational varieties.
\li The intersection ring $A^*(\ltSi(x,l))$ is generated by the pullbacks $i^*(X),i^*(L),i^*(F)$ of
the corresponding divisors in the ambient space. In particular,
let $[\ltSi(X,L,F)]\in H^*(\mXL\times\mPN)$ be the class of the lifted stratum in the
integer cohomology of the ambient space (considered as a homogeneous polynomial in $X,L,F$). Then:
\beq
A^*(\ltSi(x,l))=\frac{\mZ[i^*(X),i^*(L),i^*(F)]}
{i^*(X)^3,i^*(L)^3,i^*(X)^2-i^*(L)i^*(X)+i^*(L)^2,Ch(i^*(X),i^*(L),i^*(F))}
\eeq
Here $Ch(i^*(X),i^*(L),i^*(F))$ is a homogeneous polynomial defined by the identity
$Ch(X,L,1)[\ltSi(X,L,1)]=1$ in the cohomology ring.
\li the Picard group for a (lifted, closed) stratum of linear singularity is a free group of rank 3:
$Span_\mZ[i^*(X),i^*(L),i^*(F)]$  (except for the case of ordinary multiple point: $Span_\mZ[i^*(X),i^*(F)]$).

\eprop
The first part of this theorem was proved in \cite{Ker06}. The second follows from
the propositions \ref{ThmChowRingProjVectorBundle} and \ref{ThmSegreClassVsCohomClass}.
\\\\\\
From the intersection ring of the closed stratum one proceeds to that of the semi-compactification.
To find the corresponding Picard group, one should factor by the boundary divisors.
In this case we obtain:
\bprop
For a given singularity type $\mS$ and a "generic" degree of curve $d$ the $Pic(\Si_\mS)$ is torsion.
\eprop
Here by generic we mean that $d$ does not satisfy a specific Diophantine equation
(which is fixed uniquely by the singularity type). This equation has usually no integer solutions.

Therefore the knowledge of the full cohomology class enables us to study the intersection theory and solves the
whole class of enumerative
problems concerning the specified singularity type (e.g. consideration of curves with a restriction on the position
of singular point, on the tangent cone etc.). It also provides the necessary setup for enumerative questions related
to tangencies of singular curves.
\subsubsection{When is  the proper stratum affine}\label{SecIntroAffineStratum}
Having identified the (irreducible) boundary divisors $D_i$ and calculated their equivalence classes,
one can check whether the proper stratum is affine.
Note that in view of non-pure dimensionality of the boundary (sec. \ref{SecIntroCompactificationAndResolution})
the affinity is not at all obvious.
Somewhat unexpectedly, it appears that the result
depends on the degree of curves $d$ (i.e. its relation to parameters of the singularity).

We use the following standard criterion:
\bfact
Let $\Si$ be a quasi-projective variety, such that $\Si\subset\lSi\subset\mP^{n_1}\times..\times\mP^{n_k}$. Suppose
the boundary is pure dimensional, of codimension 1 (i.e. a true hypersurface) and its decomposition
into irreducible component is  $\bar{X}\setminus X=\bigcup_i D_i$.
\li If $\sum a_i D_i$ is ample on $\lSi$ for some $\{a_i\in\mN\}_i$ the $\Si$ is affine.
\li If $\Si$ is affine then for some $\{a_i\in\mN\}_i$ the divisor $\sum a_i D_i$ is nef (i.e. for
every effective one-dimensional cycle $C$ the intersection $C.\sum a_i D_i$ is non-negative).
\efact
The following observation is useful.
Suppose one component of the boundary $D_1$ is ample and big. Then by choosing a big $a_1$ one can assure
that $a_1D_1+\sum_{i\ge2}D_i$ is ample too. Therefore to prove that the stratum is affine it suffices
to find just one ample component of the boundary. Alternatively, in this case $\ltSi\setminus D_1$ is
already affine. And all further deletions of divisors preserve affinity.

On the other hand, to prove that a divisor $D$ is not nef, we can consider its maximal self-intersection: $D^{dim(D)+1}$.
It corresponds to a zero dimensional scheme, so if at least one coefficient of the expression is negative,
the initial divisor is not nef.

Applying this idea to the boundary component of a stratum, corresponding to the singularity type with
higher order of determinacy, we obtain the following criterion:
\bprop
Let $\mS$ be a linear singularity type, with the multiplicity $p$ and the order of determinacy $q$.
Then the proper stratum $\Si_\mS$ is affine for $d\ge\frac{3pq}{p+q}$.
\eprop
This criterion misses the region $q\le d<\frac{3pq}{p+q}$, for example for ordinary multiple point: $p\le d<\frac{3}{2}p$.
The bound can often be improved (by considering other components of the boundary). However, there is
almost always a small region $d\gtrsim q$ where the stratum is {\it non-affine}.  This is in
some contradiction to the natural guess (after \cite{DiazHarris88}) that a proper stratum is affine.
\subsubsection{Curves with two singular points}\label{SecIntroCurvesTwoSingPoints}
As an auxiliary problem we treat the problem of collision of two singular points and calculate the degrees of
some strata of the type $\Si_{\mS,A_1}$.

According to the general philosophy of Thom, proved by Kazarian (cf. \cite{Kaz2,Kaz3}), for a collection of
singularities $\mS_1..\mS_r$
the degree of the stratum $\tSi_{\mS_1..\mS_r}$ is expressed through the universal Thom polynomials $S_{\mS_i}$,
depending on the relative Chern classes of the ambient space and the linear system. The expression is:
\beq
deg(\Si_{\mS_1..\mS_r})=\sum_{J_1\bigsqcup..\bigsqcup J_k}S_{\mS_{J_1}}S_{\mS_{J_k}}
\eeq
Therefore our results give the specializations of the polynomials $S_{\mS_i}$, to the case of a complete linear system of
plane curves.

Of course the universal Thom polynomials cannot be restored from our answers. However the answers are important,
imposing already some numerical restrictions on the polynomials.
\subsubsection{Acknowledgements}
This work is a tail of my PhD, done under the supervision of E.Shustin, to whom I wish to express my gratitude.
The conversations and advices of (in alphabetic order) P.Aluffi, G.-M.Greuel, D.Markoushevitch, A.Nemethi,
D.Stepanov, I.Tyomkin,
W.Veys were highly important.

The work was done during my stay in Max Planck Institute f\"ur Mathematik, Bonn. Many thanks
for excellent working conditions.
\subsection{Example: curves with an ordinary multiple point}\label{SecIntroOrdinaryMultiplePoint}
The simplest example is the stratum of curves with ordinary multiple point $f=x^{p+1}_1+x^{p+1}_2$. In this case
the calculations are immediate, in particular the intersection ring was obtained in \cite{MiretValls}.

The lifted stratum is just the universal curve \cite{Ker06}:
\beq
\ltSi(x)={\{(x,f)|~f|_x^{(p)}=0\}}\subset\mP^2_x\times\mPN
\eeq
It is a globally complete intersection (all the defining equations are transverse). Therefore the lifted
 closed stratum is T-smooth.
Correspondingly
its cohomology class is (for the details and notations cf. section \ref{SecDefinitionsNotations})
\beq\label{EqCohClassOrdinaryMultiplePoint}
[\ltSi(x)]=\Big((d-p)X+F\Big)^{{p+2}\choose{2}}\in H^*(\mPN\times\mP^2_x)
\eeq
Here $F,X$ are the generators of the ring $H^*(\mPN\times\mP^2_x)$.

The coefficients of this polynomial have direct enumerative meaning, being related to intrinsic divisors:
\li The number of curves in a generic linear system (of appropriate dimension) with a point of multiplicity $p+1$
(lying anywhere in the plane)
is the coefficient of $X^2F^{{p+2}\choose{2}-2}$.
\li The number of curves in a generic linear system (of appropriate dimension) with a point of multiplicity $p+1$
lying on a fixed generic line, is the coefficient of $XF^{{p+2}\choose{2}-1}$.
\li The number of curves in a generic linear system (of appropriate dimension) with a point of multiplicity $p+1$
fixed in the plane, is the coefficient of $F^{{p+2}\choose{2}}$. It is obviously 1.

Thinking of the parameter space $\mPN$ as the projectivization of a vector space $\mPN=Proj(\hPN)$, the
lifted stratum is the projectivization of a vector bundle $\htSi\subset\mP^2\times\hPN$. Correspondingly
its intersection ring
is completely fixed by propositions \ref{ThmChowRingProjVectorBundle} and \ref{ThmSegreClassVsCohomClass}.
Namely, $A^*(\ltSi(x))=A^*(\mP^2_x)[F]\diagup\Big(F^r+F^{r-1}c_1+..+c_r\Big)$. Here  $X,F\in A^2(\ltSi(x))$ are the
pullbakcs of $X,F\in H^2(\mPN\times\mP^2_x)$ and the total Chern class $1+c_1..+c_r$
is fixed by the total Segre class $s=\Big(1+(d-p)X\Big)^{{p+2}\choose{2}}$. This illustrates the proposition
\ref{ThmIntroIntersRingLinearStrat}.
\subsubsection{Resolution of the vector bundle.}
In this case it is also easy to write explicit resolution for the vector bundle \cite{MiretValls}.
A curve belongs to the stratum iff $jet_p(f)=0$. Thus (in notations of section \ref{SecDefinitionsNotations})
$f\in S^pQ^*\otimes S^{d-p}(\hPt)^*$. So, we have just the standard Kozsul resolution:
\beq
0\to S^{d-p-1}(\hPt)^*\otimes\wedge^2Q^*\otimes S^{p-1}Q^*\stackrel{\al}{\to} S^{d-p}(\hPt)^*\otimes S^{p}Q^*
\stackrel{\be}{\to}\htSi(x)\to0
\eeq
with the maps $\al:~f_1\otimes(\xi_1\wedge\xi_2)\otimes f_2\to (f_1\xi_1)\otimes(f_2\xi_2)-(f_1\xi_2)\otimes(f_2\xi_1)$
and $\be:~f_1\otimes f_2\to (f_1f_2)$.

From the resolution the total Chern of the bundle $\htSi(x)$ is directly calculated
(cf. sec \ref{SecDefinitionsRelevantBundles}).
Finally (cf. \cite[proposition 1.2]{MiretValls})
\beq
A^*(\ltSi(x))=\frac{\mZ[X,F]}{\Big(X^3,~~F^r-{p+2\choose{2}}(d-p)XF^{r-1}
+\Big({p+2\choose{2}}+{{{p+2\choose{2}}\choose{2}}}\Big)(d-p)^2X^2F^{r-2}\Big)}
,~~~~r=N_d+1-{p+2\choose{2}}
\eeq
The relations in the intersection ring $A^*(\ltSi(x))$ are obtained most easily from the cohomology class
(\ref{EqCohClassOrdinaryMultiplePoint}):
\beq
X^2F^{r-1}=1,~~XF^r={p+2\choose{2}}(d-p),~~~F^{r+1}={{p+2\choose{2}}\choose{2}}(d-p)^2
\eeq

\subsubsection{Some extrinsic divisors/cycles} can be described very explicitly and
their classes can be written out immediately. For example, consider a
cycle in $\ltSi$, consisting of curves with ordinary multiple point, such that the $i$'th (smooth)
branch and its tangent line $l_i$ have tangency of order $k_i\ge2$. (Here $l_i$'s are 1-forms defining the lines.)

Start from the branch decomposition: $f=(l_1+..)(l_2+..)..(l_{p+1}+..)+\text{\it higher order terms}$.
The condition of prescribed order of tangency gives:
\beq
f=(l_1(1+..)+m_{k_1})(l_2(1+..)+m_{k_2})..(l_{p+1}(1+..)+m_{k_{p+1}})+m_{k_1+..+k_{p+1}+1}
\eeq
(here $m_i$ is the corresponding local ideal).

So the corresponding (lifted) cycle is defined as (cf. \cite{Ker06})
\beq
\ltSi_{k_1..k_{p+1}}(x,l_i,A_{k_i},B_{k_i-2})=\left\{\Big(\!\!\ber x,\{l_i\},\{A_{k_i}\}\\\{B_{k_i-2}\},f\eer\!\!\Big)
\Big|~~\ber f|_x^{\sum k_i-(p+1)}\sim SYM(A_{k_1}..A_{k_{p+1}}),\\
A_{k_i}(x)\sim SYM(l_i,B_{k_i-2}),~~~l_i(x)=0\eer\right\}
\eeq
here $A_{k_i},B_{k_i}$ are some auxiliary forms, all the notations are from section \ref{SecDefinitionsNotations}.
As always, the definition is a combination of standard proportionality conditions, which are mutually transverse.
Correspondingly the cohomology class is written immediately. By projecting out the auxiliary variables
($l_i,A_{k_i},B_{k_i}$) one obtains the class of the minimal lifting of the cycle: $[\ltSi_{k_1..k_{p+1}}(x)]$.
Finally, representing this class as the product $[\ltSi(x)][C]$ we get the expression $[C]$ for the cycle
in terms of pulled back generators $X,F$.
\subsubsection{Boundary divisors and affinity of the proper stratum.}
Consider now the semi-compactification $\bSi\subset\lSi$. The boundary divisor corresponds
to the degenerate multiple point. As the original point is characterized by its (non-coinciding) tangents,
the only minimal (topological) degeneration is: two tangents merge. Correspondingly, the boundary
(the stratum adjacent in codimension 1) is $\lSi_{x^{p+1}_1+x^2x^{p-1}_2+x^{p+2}_2}$, which is a hypersurface.

The boundary divisor is irreducible (since it represents an equisingular stratum).
Its class is $2pF+p(2d-3(p+1))X$ (cf. \cite[Appendix A.1.2]{Ker06}).
So, it is very ample for $2d\ge 3(p+1)$ (as a positive combination of very ample
divisors).  Thus, from the discussion in \ref{SecIntroAffineStratum}
it follows that both $\btSi$ and the proper stratum  $\Si$ are {\it affine} for $d\ge\frac{3(p+1)}{2}$.

To check whether for $p+1\le d<\frac{3(p+1)}{2}$ the divisor is nef, we consider its maximal self-intersection:
$\Big(2pF+p(2d-3(p+1))X\Big)^N$. Here $N$ is the maximal number such that the product is not zero.
So, the resulting variety of the intersection is zero dimensional. Thus one negative
coefficient of the product will cause non-nefness.

Direct check shows that for $d=p+1$ the coefficient of $X$ is negative for all $p$,
for $d=p+2$ this coefficient is negative if $p\ge8$ etc..
So, there exists a region near $(p+1)$ for which the semi-compactification is not affine.
\\
The Picard group of the semi-compactification is $\frac{\mZ[X,F]}{2pF+p\Big(2d-3(p+1)\Big)X}$.

Now move to the proper stratum. The boundary divisor of $\tSi\subset\btSi$ is the stratum
$\bSi_{x^{p+1}_1+x^{p+1}_2,A_1}$, i.e. the curve acquires an additional node.

Its divisor class is obtained in section \ref{SecMultiSingularities}, it is
\beq
\al X+\be F,~~\al=-pd(4+3p)+3p(2+3p+p^2),~~~~~\be=3(d-1)^2-4p-3p^2
\eeq
Correspondingly the Picard group of the proper stratum is $\frac{\mZ[X,F]}{2pF+p(2d-3(p+1))X,~\al X+\be F}$
recovering the known results of \cite{MiretXamboDescamp94} and \cite[section 5]{MiretValls}
(the later up to a misprint).

Again, this class is not nef for $d\ge p+2$.  Therefore for $p\ge8$, there exists a region
of values of $d$ (lying inside $(p+1,\frac{3(p+1)}{2})$) for which the proper stratum is {\it not affine}.

\section{Some definitions, notations}\label{SecDefinitionsNotations}
\subsection{The relevant notions and results}
\subsubsection{The ambient space and the cohomology.}
In this paper we deal with many rational equivalence classes of various varieties, embedded into various (products of)
projective spaces.
To simplify the formulae we adopt the following notation. If we denote a point in the space
$\mP^2_x$ by the letter $x$, then the homogeneous coordinates are $(x_0,x_1,x_2)$. The generator of the cohomology
(or intersection)
ring of this $\mP^2_x$ is denoted by the upper-case letter $X$, so
that $H^*(\mP^2_x)=\mZ[X]/(X^{3})$. Alternatively $X$ is the first Chern class of the dual tautological bundle
$\cO_{\mP^n}(1)$.

 By the same letter we also
denote the hyperplane class in homology of $\mP^n_x$. Since it is always clear,
where we speak about coordinates and where about (co)homology classes, no confusion arises.
To demonstrate this, consider the hypersurface
\beq
V=\{(x,y,f)|~f(x,y)=0\}\subset\mP^n_{x}\times\mP^n_{y}\times\mPN
\eeq
Here $f$ is a bi-homogeneous polynomial of bi-degree $d_x,d_y$ in homogeneous coordinates \mbox{$(x_0,\dots,x_n)$},
$(y_0,\dots,y_n)$,
the coefficients of $f$ are the homogeneous coordinates in the parameter space $\mPN$.
The cohomology class of this hypersurface is
\beq\label{DemostratCohomClasses}
[V]=d_xX+d_yY+F\in H^2(\mP^n_{x}\times\mP^n_{y}\times\mPN)
\eeq
A (projective)
line through the point $x\in\mP^2_x$ is defined by a 1-form $l$ (so that $l\in(\mP^2_l)^*,~l(x)=0)$. Correspondingly
the generator of $H^*((\mP^2_l)^*)$ is denoted by $L$.

For projective space homology and rational equivalence coincide. For other varieties we consider
rational equivalence of cycles. For a subvariety of multi-projective space  $\Si\into\mP^{n_1}\times..\times\mP^{n_r}$
we are interested in its intersection ring $A^*({\Si})$. A divisor $D$ in the ambient space gives the pullback
$[i^*(D)]\in A^*(\Si)$. To avoid messy notation we often omit the pullback sings (in this case
we specify the relevant intersection ring).
\subsubsection{The cohomology class of the diagonal}\label{SecClassDiagonal}
We often use the formula for the cohomology class of the diagonal $\Delta=\{x=y\}\subset\mP_x^N\times\mP_y^N$:
\beq
[\Delta]=\sum^N_{i=0} X^{N-i}Y^i
\eeq
For example, a condition of proportionality of two symmetric forms $f^{(p)}\sim g^{(p)}$ is just the coincidence
of the corresponding points in projective space.
\subsubsection{The cohomology classes of some degenerations}\label{SecClassDegenerations}
\parbox{15cm}
{$\bullet$ Removing a monomial. We often need to degenerate by demanding that a monomial on the Newton diagram is absent.
The condition is: a monomial $x^p_1x^q_2$ should be absent in the normal form (i.e. its coefficient must vanish).
The class of this condition was calculated in \cite[section A.1.2]{Ker06} and the corresponding degeneration is:
\beq
[\ltSi_1(x,l)]\Big(F+(d-p-2q)X+(q-p)L\Big)=[\ltSi_2(x,l)]
\eeq
}
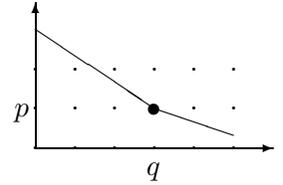
\begin{picture}(0,0)(-20,-70)
\mesh{0}{-90}{5}{2}{15}{90}{55}
\put(0,-45){\line(3,-2){45}}  \put(45,-75){\line(3,-1){30}} \put(42,-78){$\bullet$}
\put(-8,-78){$p$}  \put(42,-100){$q$}
\end{picture}
\\
This gives the class of the divisor in $\mXL\times\mPN$. The class of its pull-back to the
stratum $\ltSi_1(x,l)$ is obtained by pulling back the hyperplanes $X,L,F$, so it is:
$i^*(F)+(d-p-2q)i^*(X)+(q-p)i^*(L)$.
\\\\
\parbox{16cm}{$\bullet$ Another type of degeneration is when a quasi-homogenous form corresponding to a face of
the Newton diagram becomes degenerate. The simplest type of degeneration is the double root.
Let $a_0x^{ps}_1+..+a_sx^{qs}_2$ be the form, with $(p,q)=1$. The class of the
corresponding degeneration condition is computed by recursion:}
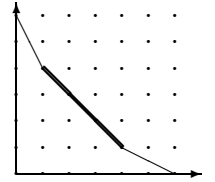
\begin{picture}(0,0)(-10,30)
\mesh{0}{0}{6}{6}{10}{70}{65}
\put(0,60){\line(1,-2){10}} \put(60,0){\line(-2,1){20}}  \thicklines \put(10,40){\line(1,-1){30}}\put(10.5,40.5){\line(1,-1){30}}
\end{picture}
\bprop
~$[degeneration]=2(s-1)\Big(i^*(F)+i^*(X)(d-s(q+\frac{p}{2}))+i^*(L)\frac{s}{2}(q-p)\Big)$
\eprop
\bpr
Let $\Si$ be the initial stratum and $D$ be the needed divisor.
Apply the degeneration procedure (as in \cite{Ker06})
to arrive at a stratum with the quasi-homogeneous form being maximally degenerated (i.e.
$x^{ps}_1$). For this demand that its coefficients vanish (one-by-one).

Let $(j,s-j)$ denote the divisor along which the coefficient of $x^{pj}_1x^{p(s-j)}_2$ vanishes.
The first degeneration is: $(\Si\cap D)\cap (0,s)=2\Big(\Si\cap(0,s)\cap(1,s-1)\Big)\cup\Big(\Si_{(0,s)}D\Big)$.
Applying the degeneration $(1,s-1)$ we get:
$\Big((\Si\cap D)\cap (0,s)-2\Big(\Si\cap(0,s)\cap(1,s-1)\Big)\Big)\cap(1,s-1)=
2\Big(\Si\cap(0,s)\cap(1,s-1)\cap(2,s-2)\Big)\cup\Big(\Si_{(0,s),(1,s-1)}D\Big)$
Continuing in this way we reach the stratum with maximally degenerate quasi-homogenous form: $x^{ps}_1$.
This stratum can be represented as $\Si\bigcap^{s-1}_{i=1}(i,s-i)$. This provides
the needed equation for the cohomology class. As all
the degenerations are invertible we divide by their classes and get:
\beq\ber
[D]=2[(1,s-1)]+....+2[(s-1,1)]=2\sum^{s-1}_{i=1}\Big(F+(d-pi-2q(s-i))X+(q(s-i)-pi)L\Big)=\\=
2(s-1)\Big(F+X(d-s(q+\frac{p}{2}))+L\frac{s}{2}(q-p)\Big)
\eer\eeq
\epr
\subsubsection{Symmetric forms.}
We often work with symmetric $p-$forms $\Omega^{p}\!\!\in\!\! S^p(\hPt)^*$ (here $(\hPt)^*$ is a 3-dimensional vector
space of linear forms).
Thinking of the form as of a symmetric tensor with $p$ indices ($\Omega^{(p)}_{i_1,\dots,i_p}$), we often write
$\Omega^{(p)}(\underbrace{x,\dots,x}_{k})$ as a shorthand for the tensor, multiplied $k$ times by the point $x\in\hPt$
\beq
\Omega^{(p)}(\underbrace{x,\dots,x}_{k}):=\sum_{0\le i_1,\dots,i_k\le2}\Omega^{(p)}_{i_1,\dots,i_p}x_{i_1}\dots x_{i_k}
\eeq
So, for example, the expression $\Omega^{(p)}(x)$ is a $(p-1)-$form. Unless stated otherwise, we assume the
symmetric form $\Omega^{(p)}$ to be generic (in particular non-degenerate, i.e. the corresponding hypersurface
$\{\Omega^{(p)}(\underbrace{x,\dots,x}_{p})=0\}\subset\mP_x^n$ is smooth).

Symmetric forms typically occur as tensors of derivatives of order $p$, e.g. $f^{(p)}$. Sometimes, to emphasize
the point at which the derivatives are calculated we assign it. So, e.g.
$f|_x^{(p)}(\underbrace{y,\dots,y}_k)$ means: the tensor of derivatives of
order $p$, calculated at the point $x$, and contracted $k$ times with $y$. Usually it will be evident, at which point
the derivative is calculated, in such cases we will omit the subscript $x$.
\subsubsection{Relevant bundles.}\label{SecDefinitionsRelevantBundles}
We constantly work with projectivization of vector bundles (vector spaces). To avoid messy notation we adopt the
following convention: if a projective bundle/space is denoted by $Z$ then the corresponding vector bundle/space
is $\hat{Z}$.
For example $\mP^n=$Proj$(\hPn)$. An equisingular stratum $\Si$ is the projectivization of the corresponding bundle
$\hSi\subset\hPN$.

The following bundles will constantly occur.
The tautological and quotient bundles are related by $0\to\cO_{\mP^n}(-1)\to\hPn|_{\mP^n}\ra Q\to0$. Here
$\cO_{\mP^n}(-1)=\{(x,v)\in\mP^n\times\hPn|v\in x\}$. The dual sequence
$0\ra Q^*\to\hPn^*|_{\mP^n}\to\cO_{\mP^n}(1)\ra 0$
defines the bundle $Q^*$ whose fiber over $x$ consists of one-forms vanishing at $x$. It is of rank $n$,
its Chern class is $c(Q^*)=\frac{1}{1+X}$.

The symmetric power $S^pQ^*$ has as its fiber polynomials (of degree $p$) vanishing at a given point up
to the order $p-1$. Its total Chern class is computed by the general rule \cite[section 3.2]{Ful}. In the case
of plane curves $n=2$ it is especially simple: $c(S^pQ^*)=1-{p+1\choose{2}}X+{{p+1\choose{2}}+1\choose{2}}X^2$

For further reference we mention also the tensor product $E\otimes\cL$ of rank $r$ bundle with a line bundle.
Its total Chern class is
\beq
c(E\otimes\cL)=1+c_1(E)+rc_1(\cL)+\Big(c_2(E)+(r-1)c_1(E)c_1(\cL)+{r\choose{2}}c^2_1(\cL)\Big)+..
\eeq

To calculate the intersection ring of the projectivization of a bundle we use
the following classical theorem:\\
\parbox{15cm}
{\bthe\label{ThmChowRingProjVectorBundle}
Let a vector bundle $E_B$ of rank $r$ and its projectivization $\mP E$ be as on the diagram,
 $\xi=c_1(\cO(1)_{\mP E})$ be the first Chern class of the dual tautological bundle.
Then $A^*(\mP E)=A^*(B)[\xi]\diagup (\xi^r+\pi^*(c_1(E))\xi^{r-1}+..\pi^*(c_r(E)))$
\ethe}
{\xymatrix{E\ar[d]\ar[dr]^{proj} &\save-<0.1cm,0cm>*{\pi^*(E)\!\!\supset}\restore\ar[d]
&\save-<0.2cm,0cm>*{\cO(-1)}\ar[dl]\restore\\B&\ar[l]_{\pi}{\mP E}}
}
\\
In our case, $\Si\into\mPN$ the first Chern class is:
$\xi=c_1(\cO_{\hSi}(1))=c_1(i^*\cO_{\hPN}(1))=i^*c_1(\cO_{\hPN}(1))=i^*(F)$.
\\
\parbox{15cm}
{To apply this theorem we need the total Chern class of the bundle $E\ra B$. In this paper we constantly
meet the situation where the total space of the bundle is embedded (cf. the diagram), so that
the fibres of $E$ are linear subspaces of $\hPN$. And
the cohomology class $[Proj(E)]\in H^*(B\times\mPN)$ is known \cite{Ker06}.
}
{\xymatrix{E\ar[d] &\save-<0.1cm,0cm>*{\subset B\times\hPN}\ar[dl]\restore\\ B
}

 This fixes the total Chern class completely:
\bprop\label{ThmSegreClassVsCohomClass}
Let the class $[Proj(E)]$ be given by a (homogeneous) polynomial $P(...,F)$ in the generators
of the intersection ring of $B\times\mPN$. Then the polynomial is monic in $F$ and the total Segre class of
the bundle $E\ra B$ is $P(..,1)$.
\eprop
\bpr
\li Let $[pt]\in H^*(B)$ be the class of a point. Then $[pt][Proj(E)]\in H^*(B\times\mPN)$ is the class dual to
the fiber over a point. So it is non-zero, being the class of a hyperplane in $\mPN$. Thus
the polynomial is monic in $F$.
\li From the previous proposition we get the identity: $F^r+\pi^*(c_1(E))F^{r-1}+..\pi^*(c_r(E))=0\in A^*(Proj(E))$.
Here $F$ is the pullback of the hyperplane class in $\mP^n_f$ while $\pi^*(c_i(E))$ are some classes on $B$.
By pushing-forward the identity
to the ambient space $B\times\mPN$ we get
$\Big(F^r+\pi^*(c_1(E))F^{r-1}+..\pi^*(c_r(E))\Big)[Proj(E)]=0\in A^{(N_d+1)}(B\times\mPN)$.

As the class $[Proj(E)]$ is monic in $F$, the Chern classes are restored uniquely from this identity.
In fact the identity can be lifted to an identity in the bigger polynomial ring :
\beq
\Big(F^r+\pi^*(c_1(E))F^{r-1}+..\pi^*(c_r(E))\Big)[Proj(E)]=F^{N_d+1}\in A^{(N_d+1)}(B)[F]
\eeq
Thus the statement for the total Segre class follows.
\epr
\subsubsection{Blowup along the diagonal.}\label{SecBlowupOverDiagonal}
Here we consider the construction and properties of the blown up product $\tpp$. It is used for the construction
of the strata $\Si_{\mS,A_1}$.

Represent it
as the incidence variety of pairs of points and lines through them:
\beq
\tpp=\{(x,y,l)|~l(x)=0=l(y)\}\into\mP^2_x\times\mP^2_y\times(\mP^2_l)^*
\eeq
The variety is a complete intersection, thus its cohomology class is the product
$[\tpp]=(L+X)(L+Y)\in H^4(\mP^2_x\times\mP^2_y\times(\mP^2_l)^*)$.

We will often need the cohomology class of the exceptional divisor $E_\De=\{x=y,~l(x)=0\}\subset\tpp$. Again, as it
is a transversal intersection of the two conditions we get:
$[E_\De]=(L+X)(X^2+XY+Y^2)\in H^6(\mP^2_x\times\mP^2_y\times(\mP^2_l)^*)$.

We also need its other form: as a class inside the blow-up $[E_\De]\in A^*(\tpp)$. Recall the intersection ring
of the blow-up. Represent it as a projectivization $\tpp=Proj(E)$, where $E$ is rank-2 (tautological) bundle over
the incidence $\{(x,l)|x\in l\}$. The incidence itself is a projectivization of the tautological
bundle over the grassmanian of lines in $\mP^2$. In total, the intersection ring of $\tpp$ is generated by the
pull-backs of the generators of that of the ambient space $\mP^2_x\times\mP^2_y\times(\mP^2_l)^*$, i.e. by
$i^*(X)$, $i^*(Y)$, $i^*(L)$. The relations arise from the relations in the cohomology of ambient space (which are:
$X^3=0=Y^3=L^3$) and the defining equations. We abuse notations and omit the pullback sings $i^*$ for compactness.

So the intersection ring is $A^*(\tpp)=\mZ[X,Y,L]/(X^3,Y^3,X^2-XL+L^2,Y^2-YL+L^2)$.

The class of the exceptional divisor in this ring is $X+Y-L$. It can be obtained, for example, by noticing
that the hypersurface $\bpm x_0&x_1\\ y_0&y_1\epm=0$ contains the exceptional divisor and also the piece $l_2=0$.

The identity $i_*(X+Y-L)=(L+X)(X^2+XY+Y^2)\in H^6(\mP^2_x\times\mP^2_y\times(\mP^2_l)^*)$ is directly verified.
\subsection{On the spaces of approximating k-jets}\label{SecSpacesOfJets}
Usually the defining conditions of singular germs are formulated in terms of k-jets of functions (where $k=o.d.-1$,
here $o.d.$ being the order of determinacy of the singularity). Therefore we consider the parameter spaces
 of the projectivized k-jets.

For a variety $Z=Spec(A)$ the classical k-jets scheme is defined as $\cL_k(Z)=Spec(A[\ep]/\ep^{k+1})$.
Correspondingly, its elements are jets of smooth or uni-branch curves. As the closure of every equisingular
stratum contains points corresponding to multi-branch singularities, we have to consider more general
spaces of jets: jets of functions. So, to a germ $(f,0)$ we assign its k-jet $b_k=jet_kf(x_1,x_2)$.

Here we consider $b_k$ both as a polynomial (in local coordinates)
and as a symmetric $k-$form (in homogeneous coordinates).
Correspondingly the definition can be written as $f|_x^{(k)}\sim b_k$.
Such a jet defines a plane curve $\{b_k(x..x)=0\}\subset\mP^2$, abusing notations we denote
it by the same letter $b_k$.

Thus, naively we define the parameter space of k-jets as the incidence variety (the basic point, the jet):
$\{(x,b_k)|~~b_k(x..x)=0\}\subset\mP^2_x\times\mP_{b_k}$.
\bex
In the simplest case ($k=1$) the curve $b_1$ is just the tangent line and the space of projectivized 1-jets is
$\mP(T^*\mP^2)(-2)=\{(x,l)|l(x)=0\}\subset\mXL$.
Projectivization turns the space of one-forms vanishing at $x$ into the lines through the point.
For higher $k$ we have a conic osculating the curve $f=0$, the cubic etc.
\eex

For a generic point of the variety of k-jets, by the natural reduction $b_k\ra b_{k'}=jet_{k'}(b_k)$ we get all
the lower jets (approximations to the curve). However for jets with singularities at $x$ this is not the case.

To carry all the information we therefore have to blowup the parameter space along the loci of jets
with an ordinary multiple point (i.e. $jet_{k'}(b_k)=0$ or $b_k(\underbrace{x..x}_{k-k'})=0$). Thus we define
\bed The parameter space of projectivized jets is
$$\mPJ_k:=\{(x,l,b_2..b_k)|l(x)=0,~~b_2(x)\sim l,..,b_k(x)\sim b_{k-1}\}
\subset\mXL\times(\mP_{b_2})^*\dots(\mP_{b_k})^*$$
\eed
Here the conditions are that the point $x$ lies on the line $l$, the line $l$ is tangent to the conic $b_2$ at $x$ etc..
\bprop\li The parameter space is smooth and the map $\mPJ_k\to\mPJ_{k-1}$ is a (projective) bundle of rank $k+1$.
\li
The map $(x,l,b_2..b_k)\to(x,l,b_2..b_{i-1},b_{i+1}..b_k)$ is the blowup over the locus $b_{i+1}(x)=0$.
\eprop
\bpr
\li First, note that each fiber is evidently a projective space. Thus one only has to show a
trivialization over some open sets and to check that the transition maps are linear on fibers.

The automorphisms of the plane ($PGL(3)$)
act linearly on the parameter space and all the orbits are isomorphic.
Therefore it is sufficient to prove local triviality for the restriction of the bundle to a
cycle in the parameter space: fixed point $x=(1,0,0)$. For a fixed $x$, consider $b_{k-1}$
as a point in a big projective space. In fact, due to the condition $b_{k-1}(x..x)=0$ the point lies
in some hyperplane. Cover this hyperplane by standard affine sets $U_{i_1..i_{k-1}}=\{b_{i_1..i_{k-1}}=1\}$.
Thinking of each such set as being $\mC^n$, there is a unique translation from any point to the origin:
$Tr_{b_{k-1}}$. This translation induces a {\it linear} transformation on the fiber $Tr^*_{b_{k-1}}(b_k)$,
which is the needed local trivialization.

Note that the trivialization is achieved by linear transformations, therefore the transition functions are also linear.
\li The projection $(x,l,b_2,..b_k)\to(x,l,b_2..b_{i-1},b_{i+1}..b_k)$ is generically 1:1 (as $b_i\sim b_{i+1}(x)$).
It is not 1:1 over the locus with $b_{i+1}(x)=0$. Here the only restriction is $b_{i}(x)\sim b_{i-1}$,
this gives the dimension of the fiber: ${i+2\choose{2}}-{i+1\choose{2}}$.
\epr
\beR
While on the total space only $PGL(3)$ acts, the group acting on a fiber is bigger. It is a subgroup of $PGL(S^k\hPt)$
that preserves the hyperplane $b_k(x..x)=0$.
\eeR
Note that the variety $\mPJ_k$ is defined as a subvariety of multi-projective space by a combination
of standard set of conditions. Therefore its cohomology class in $H^*(\mP^2_x\times\dots(\mP_{b_k})^*,\mZ)$
is written immediately:
\beq
[\mPJ_k]=[l(x)=0][b_2(x)\sim l]\dots[b_k(x)\sim b_{k-1}]
\eeq
where the classes in square brackets are just the classes of diagonals (cf. section \ref{SecClassDiagonal}).

Alternatively, from the explicit definition above we can calculate the Chow rings of the projectivized jet spaces.
For this, we represent $\hPJ_k$ as a vector bundle over $\mPJ_{k-1}$ (of rank $r=k+1$). The above
definition immediately leads to the free resolution:
\beq
0\to \cO(-1)_{\mPJ_{k-1}}\otimes Q^*\stackrel{\al}{\to} S^kQ^*\oplus\Big(\cO(-1)_{\mPJ_{k-1}}\otimes\hPt^*\Big)
\stackrel{\be}{\to}\hPJ_k\to0
\eeq
Here the maps are: $\al:~\xi\otimes q\ra\xi q\oplus(-\xi\otimes q)$ and $\be:~q\oplus\xi\otimes f\ra q+\xi f$.

Now, the total Chern class is
\beq
c(\hPJ_k)=\frac{c(S^kQ^*)c\Big(\cO(-1)_{\mPJ_{k-1}}\otimes\hPt^*\Big)}{c(\cO(-1)_{\mPJ_{k-1}}\otimes Q^*)}
\eeq
and thus is expressed in terms of the Chern classes of $Q^*$ and $\hPJ_{k-1}$. Once the total Chern class
is known, the Chow ring is determined by the theorem \ref{ThmChowRingProjVectorBundle}.

\bex For $k=1$ we have $c(\mPJ_1)=1-X+X^2$ and the intersection ring is
$A^*(\mP T^*\mP^2)=A^*(\mP^2)\diagup(L^2-LX+X^2)$.  Here $X,L$ are the pullbacks of the
corresponding first Chern classes from the ambient space $\mP T^*\mP^2\subset\mXL$
\eex
This computes the Chow rings of the projectivized jet spaces.
\subsection{On the singularity types}\label{SecSingularityTypes}
\bed
Let $(C_x,x)\subset(\mC^2_x,x)$ and $(C_y,y)\subset(\mC^2_y,y)$ be two germs of isolated curve singularities.
They are  topologically (analytically) equivalent if there exist a homeomorphism (local analytic map)
 $(\mC^2_x,x)\mapsto(\mC^2_y,y)$ mapping $(C_x,x)$ to $(C_y,y)$. The corresponding equivalence class is
 called topological (analytical) singularity type. The variety of points (in the parameter space
$\mPN$), corresponding to curves with singularity of the same (topological/analytical) type is called
the {\bf equisingular stratum}
\eed
In the case of curves the topological type can be defined also by a (simple polynomial) representative of the type:
the {\bf normal} form. For example for several simplest types
(all the notations are from \cite{AGV}, we ignore the moduli of analytic types):
\beq\scriptstyle\ber
A_k:x^2_2+x^{k+1}_1,~~D_k:x^2_2x_1+x^{k-1}_1,~~E_{6k}:x^3_2+x^{3k+1}_1,~~E_{6k+1}:x^3_2+x_2x^{2k+1}_1,~~
E_{6k+2}:x^3_2+x^{3k+2}_1\\
J_{k\ge1,i\ge0}:x^3_2+x^2_2x^k_1+x^{3k+i}_1,~~Z_{6k-1}:x^3_2x_1+x^{3k-1}_1,~~
Z_{6k}:x^3_2x_1+x_2x^{2k}_1,~~Z_{6k+1}:x^3_2x_1+x^{3k}_1\\
X_{k\ge1,i\ge0}:x^4_2+x^3_2x^k_1+x^2_2x^{2k}_1+x^{4k+i}_1,~~W_{12k}:x^4_2+x^{4k+1}_1,~~W_{12k+1}:x^4_2+x_2x^{3k+1}_1
\eer\eeq
Using the normal form $f=\sum a_{\bf I}{\bf x}^{\bf I}$ one can draw
the Newton diagram of the singularity. Namely, one marks the points ${\bf I}$ corresponding to non-vanishing
monomials in $f$, and takes the convex hull of the sets ${\bf I}+\mR_+^2$. The envelope of the convex
hull (a chain of segments-faces) is the Newton diagram.

To each singularity type we associate the normal form and its diagram. This diagram is
called {\it the Newton diagram} of the type.
\bed
The singularity type is (generalized) Newton-non-degenerate if the truncation of
its normal form polynomial to every face of the Newton polygon is non-degenerate\footnote{the truncated polynomial has no
singular points in the torus $(C^*)^2$} or the normal form can be brought to such a form by a
locally analytic transformation.
\eed
So, for \gNN singularities the types can be specified by giving their diagrams.
As we consider the topological types/equivalence, one could expect that to bring a germ to the Newton diagram of the
normal form, one needs local homeomorphisms. However for curves the locally analytic transformation always suffice.
In this paper we restrict consideration further to the types for which only linear transformations suffice.
\bed
A \gNN singular germ is called linear if it can be brought to the Newton diagram of its type by projective
transformations only (or linear transformations in the local coordinate system centered at the
singular point).
A linear stratum is the equisingular stratum, whose open dense part consists of linear germs.
The topological type is called linear if the corresponding stratum is linear.
\eed
The linear types happen to be abundant due to the following observation (\cite[section 3.1]{Ker06})
\bprop
The topological type is linear iff every segment of the Newton diagram
 has the bounded slope: $\frac{1}{2}\leq\mbox{tg}(\alpha)\leq2$.
\eprop
\bex
The simplest class of examples of linear singularities is defined by the series: $f=x^p+y^q,~~p\leq q\leq2p$.
In general, for a given series only for a few types of singularities the strata can be linear.
In the low modality cases the linear types are:
\li{Simple singularities (no moduli)}: $A_{1\le k\le3},~~D_{4\le k\le6},~~E_{6\le k\le8}$
\li{Unimodal singularities}: $X_9(=X_{1,0}),~~J_{10}(=J_{2,0}),~~Z_{11\le k\le13},~~W_{12\le k\le 13}$
\li{Bimodal}: $Z_{1,0},~~W_{1,0},~~W_{1,1},~~W_{17},~~W_{18}$
\eex
Most singularity types are nonlinear. For example if a curve has an $A_4$ point, the best we can do by projective
transformations is to bring it to the Newton diagram of $A_3$:
\beq
a_{0,2}x^2_2+a_{2,1}x_2x^2_1+a_{4,0}x^4_1
\eeq
This quasi-homogeneous form is degenerated ($a^2_{2,1}=4a_{0,2}a_{4,0}$) and by quadratic (nonlinear!) change
 of coordinates the normal form of $A_4$ is achieved.
\section{Explicit construction of some compactified strata}
\subsection{Liftings}
A natural approach to the resolution of a compactified stratum is to lift the stratum to a bigger ambient space.
As a bonus the lifted
strata often appear to be fibrations over some simple subvarieties of the jet space $\mPJ$. This enables
us to study their geometry.
For a linear singularity type $\mS$ with order of determinacy $o.d.=k+1$ the lifting is defined as:
\beq
\ltSi_\mS(x,l,b_2..b_{k}):=\overline{\left\{(x,l,b_2..b_{k},f)\Big|\bet $f_k$ has the singularity type $\mS$
\\with the $k$'th jet $b_k$\eet\right\}}\subset\mPJ_k\times\mPN
\eeq
The lifting or non-linear singularities is more complicated and is discussed later.
\bex
\li The minimal lifting $\ltSi(x)$ assigns to the curve its singular point. The projection $\ltSi(x)\to\lSi$
is 1:1 when restricted to the unisingular curves. In the part of multi-singular reduced curves the degree of
fiber jumps, but the fiber is still zero dimensional. The dimension of the fiber can be positive over non-reduced curves.
\li The next lifting $\ltSi(x,l)$ takes into account the tangent line to one of branches (a line in the plane is defined
 by a one form $l\in(\mP^2_l)^*$). The projection $\ltSi(x,l)\to\ltSi(x)$ is not 1:1 over the curves with several branches
 (of the same topological type) or over the curves with a singular point of high multiplicity.
\eex

One general remark is worth to mention. The lifting is defined by first lifting the proper stratum and
then taking the topological closure. The explicit defining equations of the points in the closure can be difficult to write.
However (even for non-linear singularities) the so defined variety always surjects the initial one:
\bprop\label{ThmLiftedStratumSurjects}
 Let $\tSi_1,\tSi_2$ be two consecutive liftings of the proper stratum, i.e. the projection
$\tSi_2\to\tSi_1$ eliminates one variable. Then the projection of closures
$\overline\tSi_2\to\overline\tSi_1$ is well defined and surjective.
\eprop
The proof follows easily from the consideration in classical topology (since both varieties are compact).
\subsubsection{Linear strata}
Here we consider some specific strata: those of linear singularities (defined in section \ref{SecSingularityTypes}).

For linear singularity we use the action of the group $PGL(3)$ on the stratum to fix the singular point
and a tangent line. As the linear singularity is completely fixed by its Newton diagram, the fiber over the pair
(the point, the line) is a linear subspace of $\mPN$ and the whole lifted stratum is a projective
bundle over the incidence variety.

So, for linear singularities the lifted stratum $\tSi(x,l)$ is already smooth. To simplify the strata we
lift them further: blowup the parameter space along the strata of points of higher multiplicities. So
we are led to the liftings $\tSi(x,l,b_2..b_{o.d.-1})$.

The analysis in \cite[section 3]{Ker06} gives:
\bprop
For a linear singularity with order of determinacy: $o.d.$, the (closure of the) lifted stratum is
a subvariety $\ltSi(x,l,b_2..b_{o.d.-1})\subset \mPJ_{o.d.-1}\times\mPN$. The stratum is the
projectivization of a vector bundle over a subvariety of $\mPJ_{o.d.-1}$. The base space and
the total space of the fibration are smooth. The projection to the base factors
through a chain of fibrations:
$\ltSi(x,l,b_2..b_{o.d.-1})\ra\ltSi(x,l,b_2..b_{o.d.-2})\ra..\ra\ltSi(x,l)\ra\mPJ_{o.d.-1}$
\eprop
The simplest example: ordinary multiple points is considered in \ref{SecIntroOrdinaryMultiplePoint}.
As in this example, the definition of a lifted stratum can be translated into a free resolution of
the corresponding vector bundle.
\subsubsection{An example of non-linear strata: $x^p_1+x^q_2,~~q>2p$. (For $p=2:$ $A_k$ series.)}
To formulate the defining conditions, start from the normal form of the type: $f=x^p_1+x^q_2$.
As a consequence one has $jet_{q-1}(f)=x^p_1$. To obtain the covariant condition, note that having fixed
the singular point
and the tangent line, the normal form is achieved here by locally analytic transformations
$x_i\ra x_i+\sum{\bf x_I}^{\bf I}$. Moreover in the non-linear shift only the terms of the degree $\le q-p$
are relevant. Therefore the defining condition can be written as $jet_{q-1}(f)=jet_{q-1}(b^p_{q-p})$, where
 $b_{q-p}$ is a polynomial of degree $q-p$ (alternatively a germ of smooth curve at the origin). To write this
 condition in a form covariant under $PGL(3)$ we pass to (symmetric) tensors of derivatives. Namely, let  $b_{q-p}$
 to denote also a symmetric tensor of order $(q-p)$ in variables $x_0..x_2$ (homogeneous coordinates).
Define the preliminary stratum
\beq
\ltSi(x,b_{q-p})=\overline{\left\{(x,b_{q-p},f)\Big|
\ber b_{q-p}\mbox{ is smooth},~~b_{q-p}(\underbrace{x..x}_{q-p})=0,~
~f|_x^{(q-1)}\!\!\!\!\sim\!\!(\underbrace{b_{q-p}..b_{q-p}}_{p})(\!\!\!\!\underbrace{x...x}_{q(p-1)+1-p^2}\!\!\!\!)\!\!\eer
\right\}}
\subset\mP^2_x\times\mP_{b_{q-p}}\times\mPN
\eeq
Note that in this form it is not easy to understand the boundary of the stratum, in particular the locus
where the curve $b_{q-p}$ is singular. Naively imposing the condition $jet_1(b_{q-p})=0$ results in
the singularity of multiplicity $Min(2p,q-p+1)$. This is certainly wrong since many types of
lower multiplicities are adjacent and the lifted stratum surjects (cf. proposition \ref{ThmLiftedStratumSurjects}).

To simplify the variety we lift it further to the space of jets. So, we define the lifted stratum:
\beq
\ltSi(x,l,b_2..b_{q-p})=\overline{\left\{\ber(x,l,b_2..b_{q-p},f)\\b_{i}\mbox{ are smooth}\eer\Big|
~\ber f^{(q-1)}|_x\sim (\underbrace{b_{q-p}...b_{q-p}}_{p})(\underbrace{x...x}_{q(p-1)+1-p^2}),~~l(x)=0,\\
b_{q-p}(x)\sim b_{q-p-1},..b_{2}(x)\sim l\eer \right\}}
\subset\mP^2_x\times\mPJ_{q-p}\times\mPN
\eeq
The projection $(x,l,b_2..b_{q-p})\to(x,l..b_i,b_{i+2}..b_{q-p})$ is the blowup with the center $b_{i+2}(x)=0$.
The total projection $(x,l,b_2..b_{q-p})\to(x,b_{q-p})$ is thus a chain of blowups with centers at
$\Big(b_{i}(x)=0\Big)_{i\ge2}$.  The lifted stratum is the strict transform of the preliminary stratum in these blowups.

The projection $\tSi(x,l,b_2..b_{q-p})\to\tSi(x,l)$ is 1:1 outside the points corresponding to types with
multiplicity$\ge p+2$.

The lifted stratum is a projective vector bundle outside the locus $b_{q-p}(\underbrace{x..x}_{q-p-1})=0$.

Though we have represented the lifted stratum by explicit equations, its cohomology class cannot be calculated
in a straightforward manner due to the presence of residual piece of big dimension over
$b_{q-p}(\underbrace{x..x}_{q-p-1})=0$. Rather one applies degenerations (cf. \cite{Ker06}).
In this case the projection $\ltSi\ra Aux$ has no good bundle structure and we cannot obtain
much geometric information about the stratum.
\section{The geometry of linear strata}
The representation of the (closure) of the linear strata as projectivized vector bundles enables
us to calculate the intersection rings of the strata (as formulated in theorem \ref{ThmIntroIntersRingLinearStrat}).

The theorem \ref{ThmChowRingProjVectorBundle}
reduces the problem to the intersection ring of the auxiliary space and the total Chern class of the bundle.
The first is calculated in proposition \ref{SecSpacesOfJets}, the second is fixed by the proposition
 \ref{ThmSegreClassVsCohomClass}, which in our case reads:

{\it Let the cohomology class of a lifted stratum $\ltSi\subset Aux\times\mP^D_f$ be given by a homogeneous polynomial
$P(X,L,B_2..B_k,F)$ in the generators of the cohomology ring $H^*(Aux\times\mP^D_f)$. Then the total
Segre class of the vector bundle $\htSi\ra Aux$ is given by $s=P(X,L,B_2..B_k,1)$.
}

The general procedure is described in introduction. We consider some more examples in section
\ref{SecExamplesOfStrata}. First we settle the question of pure dimensionality of the boundary.
\subsection{Boundary components}\label{SecBoundaryComponents}
\parbox{15cm}
{We prove that the boundary strata are hypersurfaces (i.e. of pure codimension 1), as is stated in proposition
\ref{ThmBoundaryIsHypersurface}.

For definiteness, we fix the following convention: {\it every Newton diagram intersects the coordinate axes
in the integral points}. This can always be achieved (by adding higher order terms) without changing
the singularity type.}
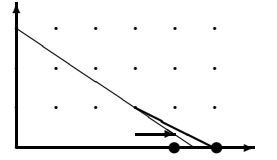
\begin{picture}(0,0)(-20,-70)
\mesh{0}{-90}{5}{3}{15}{90}{55}
\put(0,-45){\line(3,-2){67}}  \thicklines\put(45,-75){\line(2,-1){30}} \thinlines \put(45,-85){\vector(1,0){15}}
\put(57,-93){$\bullet$}   \put(73,-93){$\bullet$}
\end{picture}

Now define a partial ordering on the pairs: the Newton diagram and the stratum $(ND,\Si_{ND})$.
\bed
 We say that $(ND_1,\Si_{1})<(ND_2,\Si_{2})$ if no integral point of the diagram $ND_2$ lies under the diagram
 $ND_1$ and the topological type of $\Si_2$ is adjacent to that of $\Si_1$.

 The two pairs $(ND_1,\Si_{1})<(ND_2,\Si_{2})$ are called nearest neighbors (and the adjacency
$\Si_{S_2}\subset\lSi_{S_1}$ is called strict) if there are no intermediate pairs. For
any two pairs $(ND,\Si_{})<(ND',\Si'_{})$ the chain of nearest neighbors can be constructed:
$(ND,\Si_{ND})=(ND_1,\Si_{1})<...<(ND_k,\Si_{k})=(ND',\Si'_{})$. The chain is finite (due to the codimension, or the Milnor number).
\eed
Before we characterize possible degenerations, we consider typical examples of strict adjacency, i.e. the minimal
degenerations. We also
need the classes of the corresponding divisors in the intersection ring.
\subsubsection{Examples of degenerations}
\bex\label{ExampleDegenerationRemoveVertex}~\\
\parbox{15cm}
{Degeneration by removing a {\bf vertex} (the intersection of two faces). The condition is:
a monomial $x^py^q$, $p\ne0\ne q$ should be absent (i.e. its coefficient must vanish).
The class of this condition is given in section \ref{SecClassDegenerations} and the corresponding degeneration is:
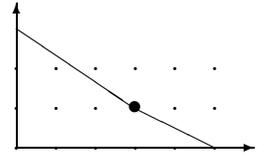
\beq
[\ltSi_1(x,l)]\Big(F+(d-p-2q)X+(q-p)L\Big)=[\ltSi_2(x,l)]
\eeq
}
\begin{picture}(0,0)(-20,-70)
\mesh{0}{-90}{5}{2}{15}{90}{55}
\put(0,-45){\line(3,-2){45}}  \put(45,-75){\line(2,-1){30}} \put(42,-77){$\bullet$}
\end{picture}
\\
This gives the class of the divisor in $\mXL\times\mPN$. The class of its pull-back to the
stratum $\ltSi_1(x,l)$ is obtained by pulling back the hyperplanes $X,L,F$, so it is:
$i^*(F)+(d-p-2q)i^*(X)+(q-p)i^*(L)$.

We immediately obtain that this boundary component is ample divisor if $d-p-2q>0$ and $q-p>0$.
\eex
\bex~\\
\parbox{15cm}
{Degeneration by removing an {\bf endpoint} (i.e. an intersection of the Newton diagram with a coordinate axis).
Let $x^r_1x^p_2+..+x^q_2$ be the quasi-homogeneous form corresponding to the face of the Newton diagram, intersecting
the $x_2$ axis.}
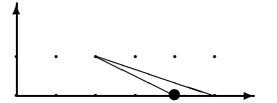
\begin{picture}(0,0)(-20,-80)
\mesh{0}{-90}{5}{1}{15}{90}{35}
\put(30,-75){\line(2,-1){30}} \put(30,-75){\line(3,-1){45}} \put(57,-92){$\bullet$}
\end{picture}
\\
 Note that if $r=1$, then this form consists of two monomials only. And the endpoint $x^q_2$
can be removed by locally analytic transformation $x_1\ra x_1+x^{q-p}_2$ (without changing the singularity type).
In such case we call the endpoint {\bf inessential}.
\eex
\bex\label{ExampleDegenerationRemoveEssentialLoosePoint}
Consider the case of {\bf essential loose endpoint} with $r>1$ and $r$ does not divide $(q-p)$.
Now, erasing the endpoint $x^q_2$ changes the topological type, we get a
strictly adjacent stratum. Its divisor class is given by the same formula as above: $i^*(F)+(d-2q)i^*(X)+qi^*(L)$.
\eex
\bex\label{ExampleDegenerationRemoveEssentialNonLoosePoint}
Another case is the {\bf essential non-loose endpoint} with $r>1$ and $r$ divides $(q-p)$. In this case the endpoint
can be removed by locally analytic shift: $x_1\ra x_1+x^{\frac{q-p}{r}}_2$ (preserving the singularity type).
So, erasing of this point does not change the topological singularity type. The possible degeneracy
of quasi-homogeneous form here gives the Newton-degenerate case, treated below.
\eex
\bex\label{ExampleDegenerationND}~\\
\parbox{15.5cm}
{ The last case is the degeneracy of a quasi-homogeneous form $x_1^{*}x_2^{*}(x_1^{ps}+..+x_2^{qs}),~~(p,q)=1$,
 corresponding to a face of the Newton diagram.In this case the type $S_2$ is typically Newton-degenerate.
 The simplest degeneracy type is:
 coincidence of the two roots, this realizes the strict adjacency.
 The cohomology class of such degeneration is calculated in section \ref{SecClassDegenerations}:
}
\begin{picture}(0,0)(-10,30)
\mesh{0}{0}{6}{6}{10}{70}{65}
\put(0,60){\line(1,-2){10}} \put(60,0){\line(-2,1){20}}  \thicklines \put(10,40){\line(1,-1){30}}\put(10.5,40.5){\line(1,-1){30}}
\end{picture}
\beq
2(s-1)\Big(i^*(F)+i^*(X)(d-s(q+\frac{p}{2}))+i^*(L)\frac{s}{2}(q-p)\Big)
\eeq
\eex
\subsubsection{General description}
In all the cases above the strict adjacency is in codimension one i.e. $\Si'\subset\lSi$ is a hypersurface.
To prove the proposition \ref{ThmBoundaryIsHypersurface} we prove that those are all the possible cases.
\bprop
Let $\Si$ be a linear stratum and the adjacency $\Si'\subset\lSi$ be strict. Then one of the following cases is
realized:
\li The type of $\Si'$ is generalized Newton-non-degenerate. Then $ND(\Si')$ is obtained from $ND(\Si)$
by either erasing a vertex (the intersection of two faces) or erasing an essential endpoint.
\li The type of $\Si'$ is Newton-degenerate. Then $ND(\Si')=ND(\Si)$, while all the quasi-homogeneous
forms, except for one, are non-degenerate. The degenerate form has one double root, all the others are single.
\eprop
\bpr Note that removing only inner points of the faces does not change the diagram, and the type remains \gNN.
So this does not change the topological type.

If at least one vertex is erased, then the topological type is changed. Thus
any further degeneration will result in a non-strict adjacency $\Si'\subset\lSi$.
The same is for an essential loose endpoint.

Suppose no intersection point and no essential loose endpoint is erased.
\li If the type of $\Si'$ is \gNN then the only possibility is that an essential non-loose endpoint of $\Si$
is removed. Since this does not change the type an additional point should be removed. It is
immediately seen that this point must be the neighbor of the essential non-loose endpoint.
Then the topological type is changed and any further degeneration will result in non-strict adjacency.
This is the scenario of example \ref{ExampleDegenerationRemoveEssentialNonLoosePoint}.
\li If the type of $\Si'$ is Newton-degenerate then at least one of the quasi-homogeneous
forms corresponding to the faces has a multiple root. The double root already changes the topological type,
thus any further degeneracy will result in a non-strict adjacency. So here the example \ref{ExampleDegenerationND}
is realized.
\epr
\subsection{Examples of strata analysis}\label{SecExamplesOfStrata}
The simplest example: ordinary multiple point was treated in section \ref{SecIntroOrdinaryMultiplePoint}.
Here we consider more examples, just repeating the analysis.\\
\parbox{16cm}
{\subsubsection{The stratum of generalized cusps: $x^p_1+x^{p+1}_2$.} (For $p=2$ it is $A_2$, for $p=3$ it is $E_6$.)
The lifted variety in this case is (cf. \cite{Ker06}):
\beq
\ltSi(x,l)=\{(x,l,f)|~f|_x^{(p)}\sim\underbrace{l\otimes ..\otimes l}_{p},~~~l(x)=0\}\subset\mXL\times\mPN
\eeq
}
\begin{picture}(0,0)(-20,0)
\mesh{0}{0}{4}{3}{10}{50}{40}
\put(0,30){\line(4,-3){40}}   \put(-10,27){\tinyT p}   \put(30,-10){\tinyT p+1}
\end{picture}

It is directly seen to be the projectivization of the corresponding vector bundle over the auxiliary space ($\mPJ_1$).
Its class $[\ltSi]\in H^*(\mPJ_1\times\mPN)$ is $(l+x)\sum Q^{{p+2\choose{2}}-1-i}(pL)^i$, with $Q=(d-p)X+F$.
This fixes the total Segre class of the bundle. Thus the Chow ring is:
\beq
A^*(\ltSi(x,l))=\frac{\mZ[i^*(X),i^*(L),i^*(F)]}
{i^*(X)^3,i^*(L)^3,i^*(X)^2-i^*(L)i^*(X)+i^*(L)^2,
F^{r-3}(F+(d-p)X-pL)
\tinyA\Big(\ber F^2-{p+2\choose{2}}(d-p)FX+\\\Big({{p+2\choose{2}}\choose{2}}-{p+2\choose{2}}^2\Big)(d-p)^2X^2
\eer\Big)}
\eeq
here $r=N_d+1-{p+2\choose{2}}$.
The boundary of $\tSi\subset\ltSi$ consists of 3 irreducible divisors:
\li Ordinary multiple point of higher multiplicity $x^{p+1}_1+x^{p+1}_2$. Its class if $F+(d-p)X-pL$
\li The type $x^p_1+x_1x^p_2+x^{p+2}_2$. Its class is $F+(d-p-1)X+(p+1)(L-X)$.
\li The curve acquires an additional node: $\Si_{x^p_1+x^{p+1}_2,A_1}$. The divisor class is
$3(d-p-1)\Big((d+p-1)F-p^2X\Big)$.

From here we immediately get:
\bprop
The Picard group of the proper stratum is
$\frac{Span_\mZ(X,F)}{\Big(3(d-p-1)((d+p-1)F-p^2X)\Big),\Big((2p+1)F+(d(2p+1)-3p(p+1))X\Big)}$
\eprop
We see that the stratum is affine for $d\ge\frac{3p(p+1)}{2p+1}$
\\
\parbox{15.5cm}
{\subsubsection{The stratum of points of type: $x^p_1+x_1x^p_2+x^{p+2}_2$} (For $p=2$ it is $A_3$, for $p=3$ it is $E_7$.)
The lifted variety in this case is:
\beq
\ltSi(x,l,b_p)=\{(x,l,b_p,f)|~f|_x^{(p+1)}\sim SYM(l,b_p),~~b_p(x)\sim(\underbrace{l..l}_{p-1}),~~~l(x)=0\}\subset\mXL\times\mP_{b_p}\times\mPN
\eeq}
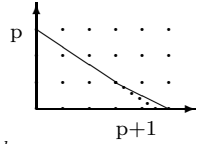
\begin{picture}(0,0)(-30,0)
\mesh{0}{0}{5}{3}{10}{60}{40}
\put(0,30){\line(3,-2){30}} \multiput(28.5,9.5)(3,-2){6}{.}  \put(30,10){\line(2,-1){20}} \put(-10,27){\tinyT p}   \put(30,-10){\tinyT p+1}
\end{picture}

(Note that here one need not consider the whole parameter space of jets $\mPJ_p$.) From this representation
the cohomology class of $\ltSi(x,l,b_p)$ is obtained immediately, and from it one gets (by projection) the class
$[\ltSi(x,l)]$.

Again, the stratum $\ltSi(x,l,b_p)$ is the projectivization of the corresponding vector bundle over the auxiliary
space ($\mPJ_1$). Thus the Chow ring is:
\beq
A^*(\ltSi(x,l))=\frac{\mZ[i^*(X),i^*(L),i^*(F)]}
{i^*(X)^3,i^*(L)^3,i^*(X)^2-i^*(L)i^*(X)+i^*(L)^2,
\frac{1}{[\ltSi(x,l)]}}
\eeq
here in the denominator by $\frac{1}{[\ltSi(x,l)]}$ we mean the polynomial expansion obtained after substitution $F\ra1$.

The boundary of $\tSi\subset\ltSi$ consists of 4 irreducible divisors:
\li The type obtained by degenerating $..x_1x^p_2\ra0$. The normal form: $x^p_1+x_1x^p_2+x^{p+2}_2$.
 Its class is $\Big(F+(d-2p-1)X+(p-1)L\Big)$.
\li The type $x^{p+1}_1+x_1x^p_2+x^{p+2}_2$, obtained by increasing multiplicity $...x^{p}_1\ra0$. Its class is
$\Big(F+(d-p)X-pL\Big)$:
\li The curve acquires an additional node: $\Si_{x^p_1+x_1x^p_2+x^{p+2}_2,A_1}$. The divisor class is
$\Big(F(3d^2-6d-1-3p^2)+(4+4p+3p^2-(d-p)(4+3p^2))X-4(p-1)L\Big)$.

From here we immediately get:
\bprop
The Picard group of the proper stratum is
$$\frac{Span_\mZ(X,F)}{\Big(3(d-p-1)(d+p-1)F-p(4-3p+3dp)X\Big),\Big((1-2p)F+((d-p)(1-2p)+p(p+1))X\Big)}$$
\eprop
\parbox{15.5cm}
{\subsubsection{The stratum of points of type: $x^{p+1}_1+x^2_1x^{p-1}_2+x^{p+2}_2$} (For $p=1$ it is $A_3$,
for $p=2$ it is $D_5$.)
The lifted variety in this case is:
\beq
\ltSi(x,l,b_p)=\{(x,l,f)|~~f|_x^{(p+1)}\sim SYM(l,l,b_p),~~b_p(x)=0,~~l(x)=0\}\subset\mXL\times\mP_{b_p}\times\mPN
\eeq}
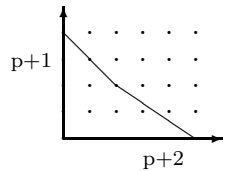
\begin{picture}(0,0)(-30,20)
\mesh{0}{0}{5}{4}{10}{60}{50}
\put(0,40){\line(1,-1){20}} \put(20,20){\line(3,-2){30}} \put(-20,27){\tinyT p+1}   \put(30,-10){\tinyT p+2}
\end{picture}

(Here $b_p$ is just an auxiliary p-form.) From this representation
the cohomology class of $\ltSi(x,l,b_p)$ is obtained immediately, and from it one gets (by projection) the class
$[\ltSi(x,l)]$.

Again, the stratum $\ltSi(x,l,b_p)$ is the projectivization of the corresponding vector bundle over the auxiliary
space ($\mPJ_1$), its Chow ring is written as before.

The boundary of $\tSi\subset\ltSi$ consists of 4 irreducible divisors:
\li The type obtained by degenerating $..x^2_1x^{p-1}_2\ra0$. The normal form: $x^{p+1}_1+x^3_1x^{p-2}_2+x^{p+2}_2$.
 Its class is $\Big(F+(d-2p)X+(p-3)L\Big)$.
\li The type obtained by degenerating $..x^{p+2}_2\ra0$. The normal form: $x^{p+1}_1+x^2_1x^{p-1}_2+x^{p+3}_2$.
 Its class is $\Big(F+(d-2(p+2))X+(p+2)L\Big)$.
 \li The type obtained by degenerating the homogeneous part $x^{p+1}_1+..+x^2_1x^{p-1}_2$.
 The resulting singularity is Newton degenerate.
 Its class is $(p-2)\Big(2F+3(d-(p-1))X\Big)$.
\li The curve acquires an additional node: $\Si_{x^{p+1}_1+x^2_1x^{p-1}_2+x^{p+2}_2,A_1}$. The divisor class is
$\Big(F(3d^2-6d-1-4p-3p^2)+(8+8p+5p^2-(d-p)(2+4p+3p^2))X-2(p-2)L\Big)$.

From here we immediately get:
\bprop
The Picard group of the proper stratum is obtained by factorization of $Span_\mZ(X,L,F)$ by the above classes.
\eprop
\section{Curves with two singular points}\label{SecMultiSingularities}
Here we consider strata of curves with two singular points, one of them being a node ($A_1$). The importance
of the strata is due to their adjacency (in codimension 1) to the strata of uni-singular curves.
The minimal lifting is now to the space of pairs of points and lines through them: $\tSi_{\mS,A_1}(x,y,l)$. In
course of the calculations we constantly face the question of collision of two singularities.
It is described by taking the flat limit.
\subsection{On the collision of singular points}\label{SecCollisionOfSingularPoints}
Working with the strata one often  faces the following problem: suppose a stratum of multi-singular
curves is defined by an explicit system of equations away from diagonals (i.e. loci where two or more singularities
merge). Describe the diagonal loci (i.e. obtain the corresponding defining equations).

As we consider the lifted strata, we are interested in the situation when one singular point is fixed.
In addition we always fix the line joining the singular points, so the collision is always along a fixed line.

The problem amounts to taking the flat limit of family and is treated by standard methods (Gr\"obner bases).
Here we consider several typical examples of merging two singular points of a plane curve. As always the
degree of curves is assumed to be sufficiently high.

We start from the product of two projective planes blown up along the diagonal:
$\tpp=\{(x,y,l)|~l(x)=0=l(y)\}\subset\mP^2_x\times\mP^2_y\times(\mP^2_l)^*$. We
usually need also the local coordinates on the exceptional divisor. In addition to the line $l$
we take also a vector $v$ through the point $x$, pointing in the $l$'th direction (i.e. $l(v)=0$).

Suppose, outside the diagonal the stratum is defined by a system of equations $\{f_\al|_x\}_\al$ (in $x$) and
 $\{g_\be|_y\}_\be$ (in $y$).
Near the diagonal expand $y=x+\ep v$ (where $\ep$ is small) and expand the equations:
\beq
g_\be|_y\ra g_\be|_x+\ep g^{(1)}_\be|_x(v)+\ep^2 g^{(2)}_\be|_x(v,v)+...
\eeq
We can assume that each expansion starts from a generically non-zero term (otherwise omit it and
divide the whole series by $\ep$). Consider now the module of syzygies of the first terms in the series,
i.e. the nontrivial relations $\sum_\be r_\be g_\be|_x=0$. Every such relation gives a series
$0+\ep\sum r_\be g^{(1)}_\be|_x(v)+..$ which is then normalized (i.e. divided by $\ep$) and added to the ideal.

Now the syzygies of the enlarged ideal are considered, and so on. The process stops when the syzygies
do not produce any new series.

Once the full ideal is obtained, one takes the limit $\ep\to0$ (just omitting the higher order
terms in all the expansions). This gives the defining ideal of the locus over the diagonal.

We consider the simplest example.
\subsubsection{Collision of two ordinary multiple points.}
\parbox{14cm}
{Suppose the multiplicities are $p+1,q+1$ such that $p\ge q$. Note that here the answer can be obtained
from the geometry of the plane blown up at the point of multiplicity $(p+1)$, cf. the picture.
The stratum is defined by
\beq
\tSi(x,y,l)=\overline{\{(x,y,l,f)|~x\ne y,~~~l(x)=0=l(y),~~f|_x^{(p)}=0=f|_y^{(q)}\}}
\eeq
Expand the equations involving $y$ around $x$, with $v=x-y$ and a small parameter $\ep$. First several terms
in the expansion vanish $f^{(q)}|_x=0=f^{(q+1)}|_x(v)=..=f^{(p)}|_x(v..v)$. This gives the (tensor) series
}
\begin{picture}(0,0)(-30,10)
\qbezier(0,10)(0,0)(0,-10)  \qbezier(-4,10)(0,0)(4,-10) \qbezier(4,10)(0,0)(-4,-10)
 \qbezier(-10,10)(20,-27)(60,10)    \qbezier(-10,-10)(20,27)(60,-10)   \qbezier(-10,0)(0,0)(60,0)
\put(-5,-15){\tinyT p}  \put(45,-10){\tinyT q}  \put(-3,15){$\downarrow$}
\put(0,25){\line(0,1){50}}  \qbezier(-10,30)(0,30)(10,30) \qbezier(-10,70)(0,70)(10,70)
 \qbezier(-10,40)(40,40)(60,60)   \qbezier(-10,50)(40,50)(60,50)   \qbezier(-10,60)(40,60)(60,40)
  \put(-20,50){\tinyT p}  \put(45,40){\tinyT q}
 \put(70,45){$\rightarrow$}
\put(100,25){\line(0,1){50}}  \qbezier(90,30)(100,30)(110,30) \qbezier(90,70)(100,70)(110,70)
 \qbezier(90,40)(100,50)(110,60)   \qbezier(90,50)(100,50)(110,50)   \qbezier(90,60)(100,50)(110,40)
 \put(97,15){$\downarrow$}
\qbezier(90,0)(100,0)(110,0)
 \qbezier(90,-10)(100,10)(110,-10)  \qbezier(90,10)(100,-10)(110,10)
\qbezier(95,13)(100,0)(105,-13)   \qbezier(95,-13)(100,0)(105,13)
\end{picture}
\beq
f|_x^{(p+1)}(\underbrace{v..v}_{p+1-q})+\ep f|_x^{(p+2)}(\underbrace{v..v}_{p+2-q})+
\ep f|_x^{(p+3)}(\underbrace{v..v}_{p+3-q})+\dots
\eeq
The syzygies are obtained as a consequence of the Euler identity for homogeneous polynomial $\sum x_i\di_if=deg(f)f$.
By successive contraction of the tensor series with $x$ we get the series
\\
\begin{tabular}{@{}>{$}c<{$}>{$}c<{$}>{$}c<{$}>{$}c<{$}>{$}c<{$}>{$}c<{$}>{$}c<{$}}
\\
\tinyA f|_x^{(p+1)}(\underbrace{v..v}_{p+1-q}) +&\tinyA \ep f|_x^{(p+2)}(\underbrace{v..v}_{p+2-q})
+&\tinyA \ep^2 f|_x^{(p+3)}(\underbrace{v..v}_{p+3-q})+&\tinyA \ep^3 f|_x^{(p+4)}(\underbrace{v..v}_{p+4-q})+&..
\\
\tinyA(d-p-2)f|_x^{(p+1)}(\underbrace{v..v}_{p+2-q}) +& \tinyA\ep (d-p-3)f|_x^{(p+2)}(\underbrace{v..v}_{p+3-q})
+&\tinyA \ep^2(d-p-4) f|_x^{(p+3)}(\underbrace{v..v}_{p+4-q})+
&\tinyA \ep^3(d-p-5) f|_x^{(p+4)}(\underbrace{v..v}_{p+5-q})+..
\\
..&..&..&..&
\\
\tinyA\prod^{q+1}_{i=2}(d-p-i)f|_x^{(p+1)}(\underbrace{v..v}_{p+1})+&\tinyA\ep\prod^{q+1}_{i=2}(d-p-1-i) f|_x^{(p+2)}(\underbrace{v..v}_{p+2})
 +&\tinyA \ep^2 \prod^{q+1}_{i=2}(d-p-2-i) f|_x^{(p+3)}(\underbrace{v..v}_{p+3})+&
 \tinyA \ep^3 \prod^{q+1}_{i=2}(d-p-3-i) f|_x^{(p+4)}(\underbrace{v..v}_{p+4})+&..
\end{tabular}
\\\\
Here the first row is the initial series, the second is obtained by contraction with $x$, the $p+2$'th
row is obtained by contraction with $\underbrace{x\times \times x}_{p+1}$.

Apply now the Gaussian elimination, bring to the upper triangular form.
\li Eliminate from the first column all the entries of the rows $2..(p+2)$. For this contract the first
row sufficient number of times with $v$ and subtract.
\li Eliminate from the second column all the entries of the rows $3..(p+2)$.
\li ...
\\Normalize the rows (i.e. divide by the necessary power of $\ep$).

In this way we get an "upper triangular" system of series (we omit the numerical coefficients):
\\
\begin{tabular}{@{}>{$}c<{$}>{$}c<{$}>{$}c<{$}>{$}c<{$}>{$}c<{$}>{$}c<{$}>{$}c<{$}}
\\
\tinyA f|_x^{(p+1)}(\underbrace{v..v}_{p+1-q}) +&\tinyA \ep f|_x^{(p+2)}(\underbrace{v..v}_{p+2-q})
+&\tinyA \ep^2 f|_x^{(p+3)}(\underbrace{v..v}_{p+3-q})+&\tinyA \ep^3 f|_x^{(p+4)}(\underbrace{v..v}_{p+4-q})+&..
\\
0 +& \tinyA f|_x^{(p+2)}(\underbrace{v..v}_{p+3-q})+&\tinyA \ep  f|_x^{(p+3)}(\underbrace{v..v}_{p+4-q})+
&\tinyA \ep^2 f|_x^{(p+4)}(\underbrace{v..v}_{p+5-q})+..
\\
0 +& 0+&\tinyA f|_x^{(p+3)}(\underbrace{v..v}_{p+5-q})+
&\tinyA \ep f|_x^{(p+4)}(\underbrace{v..v}_{p+6-q})+..
\\
..&..&..&..&
\\
0+&0 +&0+&...&+ \tinyA  f|_x^{(p+q+1)}(\underbrace{v..v}_{p+q+1})+&..
\end{tabular}
\\
Now take the limit $\ep\to0$ (i.e. just omit the higher order terms in each row) and we get the defining
system of (local) equations:
\beq
f|_x^{(p)}=0,~~~f|_x^{(p+1)}(\underbrace{v..v}_{p+1-q})=0,~~~f|_x^{(p+2)}(\underbrace{v..v}_{p+3-q})=0,
~~f|_x^{(p+3)}(\underbrace{v..v}_{p+5-q})=0~~..,f|_x^{(p++q+1)}(\underbrace{v..v}_{p+q+1})=0
\eeq
Correspondingly, the normal form of the singularity that is obtained when the two singular points merge is
\beq
x^{p+1}_1+x^{q+1}_1x^{p-q}_2+x^{q}_1x^{p+2-q}_2+...+x^{p+q+2}_2
\eeq
In several simplest cases we have: $A_1+A_1\ra A_3$, $D_4+A_1\ra D_6$, $X_9+A_1\ra X_{1,2}$, $D_4+D_4\ra J_{10}$,
 $X_9+D_4\ra Z_{13}$.
\\\\
\beR
In the case above we assumed that the two singular points collide generically, i.e. no tangents to
the branches coincide. We could safely do this, because coincidence of some tangents imposes here additional
condition on the locus over the diagonal. The resulting special locus will be therefore of a smaller dimension
(in fact it lies in the closure of the locus of generic situation).
\eeR
\subsubsection{Some other cases}
\parbox{15cm}
{In general one can use either the approach via geometry of the blown-up plane, or the analytic one.
The geometry of blowup is useful e.g. for the case of collision of two singular points: $(p_x,r_x)$
and $(p_y,r_y)$ such that the point $x$ possesses $r_1$ smooth branches with tangents different among themselves and
also different from all other tangents, and $r_1\ge r_y$.}
\begin{picture}(0,0)(-20,-10)
\qbezier(0,20)(0,0)(0,-20)  \qbezier(-10,-16)(-5,-18)(0,-18) \qbezier(-10,-20)(-5,-18)(0,-18)
 \qbezier(-10,16)(0,20)(10,16)    \qbezier(-10,20)(0,16)(10,20)

\qbezier(-10,10)(10,10)(15,0)  \qbezier(-10,0)(8,0)(10,5)   \qbezier(10,5)(5,-10)(25,0)

\qbezier(-10,-10)(12,-2)(12,12)  \put(35,0){$\ra$}
\end{picture}
\begin{picture}(0,0)(-90,-10)
\qbezier(0,20)(0,0)(0,-20)  \qbezier(-10,-16)(-5,-18)(0,-18) \qbezier(-10,-20)(-5,-18)(0,-18)
 \qbezier(-10,16)(0,20)(10,16)    \qbezier(-10,20)(0,16)(10,20)

\qbezier(-20,10)(0,10)(5,-2)  \qbezier(-20,0)(-2,0)(0,5)   \qbezier(0,5)(-5,-10)(-7,-15)

\qbezier(-20,-10)(2,-2)(2,12)
\end{picture}

In this case, in particular the resulting
multiplicity is $p_x+p_y-r_y$.
\bex Consider briefly several cases of collisions of a node with some simplest singularity types.
There will be always a preferred tangent line to the singularity, we denote it by $\tl$. To imitate
the direction of the line we introduce another point on the line $\tv$ so that $\tl(\tv)=0$.
\li The collision of $x^p_1+x^{p+1}_2$ (generalized cusp) and a node. If the collision line coincides
with the tangent to the generalized cusp ($l=\tl$) then the resulting singularity has the normal form
$x^p_1+x^2_1x^{p-2}_2+x_1x^{p+1}_2+x^{p+3}_2$ with $l$ as the tangent.
For example, for $p=2$ one has: $A_2+A_1\ra A^{(l)}_4$.

In the generic situation (i.e. $l\ne\tl$) the normal form of the resulting singularity is
$x^{p+1}_1+x^2_1x^{p-1}_2+x^{p+2}_2$ and the tangent is $l$. For example, for $p=2$ one has: $A_2+A_1\ra D_5$.
\li The collision of $x^{p-1}_1x_2+x^{p+1}_2$ and a node. If $l=\tl$ then the resulting singularity has the normal form
$x^{p-1}_1x_2+x^2_1x^{p-1}_2+x_1x^{p+1}_2+x^{p+3}_2$ with $l$ as the tangent.
For example, for $p=3$ one has: $D_5+A_1\ra D^{(l)}_7$.

In the generic situation (i.e. $l\ne\tl$) the normal form of the resulting singularity is
$x^{p+1}_1+x_1x^{p-1}_2+x^{p+2}_2$ and the tangent is $l$. For example, for $p=2$ one has: $D_5+A_1\ra X_9$,
with one of the tangents fixed.
\li The collision of $x^{p}_1+x_1x^{p}_2+x^{p+2}_2$ and a node. If $l=\tl$ then the resulting singularity
has the normal form
$x^{p}_1+x^2_1x^{p-1}_2+x_1x^{p+2}_2+x^{p+4}_2$ with $l$ as the tangent.
For example, for $p=2$ one has: $A_3+A_1\ra A^{(l)}_7$.

In the generic situation (i.e. $l\ne\tl$) the normal form of the resulting singularity is
$x^{p}_1x_2+x^2_1x^{p-1}_2+x^{p+2}_2$ and the two tangents are $l=\{x_1=0\}$ and $\tl=\{x_2=0\}$.
For example, for $p=2$ one has: $A_3+A_1\ra D_5$.
\li The collision of $x^{p-1}_1+x_1x^{p}_2+x^{p+2}_2$ and a node. If $l=\tl$ then the resulting singularity
has the normal form
$x^{p-1}_1x_2+x^2_1x^{p-1}_2+x_1x^{p+1}_2+x^{p+3}_2$ with $l$ as the tangent.
For example, for $p=3$ one has: $D_6+A_1\ra D^{(l)}_8$.

In the generic situation (i.e. $l\ne\tl$) the normal form of the resulting singularity is
$x^{p}_1x_2+x_1x^{p}_2$ i.e. the ordinary multiple point with the two fixed tangents: $l=\{x_1=0\}$ and $\tl=\{x_2=0\}$.
For example, for $p=2$ one has: $D_6+A_1\ra X_9$.
\eex
\subsection{On the classes of strata}
Here we give some examples of enumeration of curves with a node and some other singular point. According
to the Thom philosophy (cf. section \ref{SecIntroCurvesTwoSingPoints}) the answers are represented in the form:
$deg(\Si_{\mS,A_1})=S_{\mS}S_{A_1}+S_{\mS,A_1}$. As the polynomials for one singular point are known
(e.g. \cite{Kaz4} or \cite{Ker06}, and $S_{A_1}=3(d-1)^2$), we only present the part of $S_{\mS,A_1}$.
(I.e. we present the specializations of Thom polynomials to the complete linear system of plane singular curves).
\subsubsection{$x^{p+1}_1+x^{p+1}_2,A_1$}
We start from the stratum of ordinary multiple points $\tSi(x)=\{(x,f)|~f|_x^{(p)}=0\}$. Lift it
to the blowup of two planes over the diagonal $\tpp$ (discussed in \ref{SecBlowupOverDiagonal}). Intersect
with the nodal conditions at the point $y$, one-by-one ($f|_y^{(1)}=0$).

As we work in the blown up space, we should find the strict transform of the nodal conditions. For
this expand near the diagonal $d_if|_y=0+..0+\ep^{p}d_if|_x^{(p)}(\underbrace{v..v}_{p})+..$. Here $y=x+\ep v$,
where $\ep\to0$ on the diagonal, the condition is expanded up to the first term that is generically non-zero.
Thus the class of the strict transform is
\beq
\tilde{[\di_if|_y]}=[\di_if|_y]-p[E]
\eeq
Note, that the first non-zero term in the expansion above
actually vanish over some codimension 2 subvariety of the ambient space (i.e. if we consider
$\di_if|_x^{(p)}(\underbrace{v..v}{p})=0$ as an equation for $v$ not for $f$). Alternatively, the
strict transforms are not transversal there.

Correspondingly, after the intersection with the strict transform, the projection $\ltSi(x,y,l)\to\tpp$
is no longer a locally trivial fibration. The fiber over the diagonal changes (and its dimension can jump).
So, the behavior over the diagonal should be carefully checked.

Intersect the (pull-back of the) initial stratum with the strict transforms
\beq
\tSi(x)\cap\Big((\di_0f|_y=0)-pE\Big)\cap\Big((\di_1f|_y=0)-pE\Big)\cap\Big((\di_2f|_y=0)-pE\Big)
\eeq
Consider now the piece over the diagonal. It is defined by
\beq
\tSi(x)\cap\bigcap_i\Big(\di_if|_x^{(p)}(\underbrace{v..v}_{p})=0\Big)=
\{(x,f)|~f|_x^{(p)}=0=f|_x^{(p+1)}(\underbrace{v..v}_{p})\}
\eeq
This is the stratum of a linear singularity (with the normal form $x^{p+1}_1+x^2_1x^{p-1}_2+x^{p+2}_2$).

Note that its codimension equals the codimension of the needed stratum, i.e. the variety is reducible
and the contribution of this residual part should be removed. The residual piece enters with the multiplicity 2, since
it is obtained as the intersection of 3 hypersurfaces, one of which is (simply) tangent to the
intersection of the others.

The cohomology class of the residual piece is readily obtained (using the methods of \cite{Ker06}).
The cohomology of the total piece is just the product of the classes of hypersurfaces.
Thus we get the class of the needed lifted stratum $[\tSi(x,y,l)]\in H^*()$.

The degree of the actual stratum is obtained by applying Gysin homomorphism (which in this case means
just to extract the coefficient of $X^2Y^2L^2$). We get:
\beq\ber
deg(\Si_{x^{p+1}_1+x^{p+1}_2,A_1})=9{p+3\choose{4}}(d-p)^2(d-1-p)(d-1+p)-3{p+2\choose{3}}(p^2+6p+4)(d-p)^2+
3{p+2\choose{3}}(d-p)(6+5p)\\
S_{x^{p+1}_1+x^{p+1}_2,A_1}=-{p+2\choose{3}}\frac{3(d-p)}{8}\Big((d-p)(4+3p)(4+3p+p^2)-4(6+5p)\Big)
\eer\eeq

In the case of $A_1,A_1$ the answers should be divided by 2, since such is the projection $\ltSi(x,y)\ra\Si$.
\subsubsection{Some other cases}
The general strategy for the types $\mS,A_1$ is to degenerate the type $\mS$ to a ordinary multiple point.
Alternatively one can try and reach the type $\mS$ starting to degenerate from a multiple point.
In the process of degeneration residual pieces over the diagonal constantly appear and should be removed.

The method od degeneration has been described in details in \cite{Ker06}. Recall that the cohomology
class of the degenerating divisor that forces to vanish the coefficient of $x^p_1x^q_2$ is
$(d-p-2q)X+F+(q-p)L$.

We
consider some examples, following the chain of degenerations:
\beq\ber
x^{p-1}_1x_2+x^{p+1}_2 \ra x^{p}_1+x^{p+1}_2\ra x^{p+1}_1+x^{p+1}_2\\
x^{p-1}_1x_2+x_1x^p_2+x^{p+2}_2\stackrel{p\ge3}{\to} x^{p}_1+x_1x^p_2+x^{p+2}_2 \ra x^{p+1}_1+x_1x^p_2+x^{p+2}_2~
(\bet multiple point with\\ a prescribed tangency\eet)
\eer\eeq
\li $x^{p}_1+x^{p+1}_2,A_1\ra x^{p+1}_1+x^{p+1}_2,A_1.$  (For $p=2$: $A_2,A_1$, for $p=3$: $E_6,A_1$).
Degenerate by demanding that $f|_x^{(p)}=0$. As the result get the stratum $\tSi_{x^{p+1}_1+x^{p+1}_2,A_1}$
and a residual piece over the diagonal. The piece occurs because the restriction to the diagonal has a
component of curves with a point of multiplicity $p+1$ (as discussed in \ref{SecCollisionOfSingularPoints}).

So, this piece should be subtracted. Its multiplicity is 2, since the original stratum is singular
along the diagonal. One gets:
\beq
S_{x^{p}_1+x^{p+1}_2,A_1}=-\frac{3}{8}p^3\Big(p(3+p)(d-p)^2(p^2+3p-2)+4(p-1)(d-p)(p^2+3p-2)-8p\Big)
\eeq
\\
\li $x^{p}_1+x_1x^p_2+x^{p+2}_2 \ra x^{p+1}_1+x_1x^p_2+x^{p+2}_2.$ (For $p=2$: $A_3,A_1$, for $p=3$: $E_7,A_1$)
First we should calculate the class of the (auxiliary) stratum $\tSi_{x^{p+1}_1+x_1x^p_2+x^{p+2}_2,A_1}$.
This is done by degeneration of the ordinary multiple point. At this step no residual piece over
the diagonal is produced.

Next, the type $x^{p}_1+x_1x^p_2+x^{p+2}_2$ is degenerated (by demanding that the coefficient of $x^{p}_1$ vanish).
A residual piece occurs over the diagonal (collision of $x^{p}_1+x_1x^p_2+x^{p+2}_2$ and a node), its type is:
$x^{p}_1x_2+x^2_1x^{p-1}_2+x^{p+2}_2$. Subtracting
its cohomology class (with multiplicity 2) gives:
\beq
S_{x^{p}_1+x_1x^p_2+x^{p+2}_2,A_1}=-9{p+3\choose{4}}p(d-p)^2(4+p+2p^2)-\frac{3}{2}p^2(3+p)(d-p)(p^3-3p^2-p-8)+
3p(p^4+3p^3+3p^2+4p-4)
\eeq
\li $x^{p}_1+x^{p+1}_2,A_1\ra x^{p+1}_1+x^{p+1}_2,A_1.$  (For $p=2$: $A_2,A_1$, for $p=3$: $E_6,A_1$).
Degenerate by demanding that $f|_x^{(p)}=0$. As the result get the stratum $\tSi_{x^{p+1}_1+x^{p+1}_2,A_1}$
and a residual piece over the diagonal. The piece occurs because the restriction to the diagonal has a
component of curves with a point of multiplicity $p+1$ (as discussed in \ref{SecCollisionOfSingularPoints}).

So, this piece should be subtracted. Its multiplicity is 2, since the original stratum is singular
along the diagonal. One gets:
\beq
S_{x^{p}_1+x^{p+1}_2,A_1}=-\frac{3}{8}p^3\Big(p(3+p)(d-p)^2(p^2+3p-2)+4(p-1)(d-p)(p^2+3p-2)-8p\Big)
\eeq
\\
\li $x^{p+1}_1+x^2_1x^{p-1}_2+x^{p+2}_2.$ (For $p=1$: $A_3,A_1$, for $p=2$: $D_5,A_1$.)
This stratum is treated  by degeneration of $\tSi_{x^{p+1}_1+x^{p+1}_2,A_1}$. The first degeneration
$x^{p+1}_1+x^{p+1}_2\ra x^{p+1}_1+x_1x^{p}_2+x^{p+2}_2$ brings no residual piece. The second
$x^{p+1}_1+x_1x^{p}_2+x^{p+2}_2\ra x^{p+1}_1+x^2_1x^{p-1}_2+x^{p+2}_2$ brings the residual piece
over the diagonal when the tangent to the singular point at $x$ passes through the point $y$ also. This piece is:
$x^{p+1}_1+x^2_1x^{p-1}_2+x^{p+3}_2$, as always it should be subtracted with multiplicity 2.
\beq
S_{x^{p+1}_1+x^2_1x^{p-1}_2+x^{p+2}_2,A_1}=\tinyA\ber-\frac{(d-p)^2}{4}(3p^2+p+2)(p^2+3p+6)(p^2+3p+4)(p+1)\\+
\frac{d-p}{2}(p^2+3p+4)(3p^4+23p^3+30p^2+28p+12)-12-30p-28p^2-21p^3-5p^4
\eer
\eeq
\li $x^{p+1}_1+x^2_1x^{p-1}_2+x^{p+3}_2.$ (For $p=2$: $D_6,A_1$.)
The cases $p=1,2$ are exceptional here due to specific coincidences on Newton diagram.
This stratum is treated  by degeneration of $\tSi_{x^{p+1}_1+x^2_1x^{p-1}_2+x^{p+2}_2,A_1}$.
There residual piece over the diagonal occurs when the tangent line at $x$ passes through $y$. The type is:
$x^{p+1}_1+x^2_1x^{p-1}_2+x^{p+4}_2$ therefore it is irrelevant (by codimension).
\beq
S_{x^{p+1}_1+x^2_1x^{p-1}_2+x^{p+3}_2,A_1}=\tinyA\ber-\frac{(d-p)^2}{4}(p^2+3p+6)(p^2+3p+8)(9p^3+18p^2+16p+8)\\+
3\frac{d-p}{2}(p+1)(p^2+3p+6)(3p^3+20p^2+18p+16)-3(32+72p+78p^2+51p^3+10p^4)\eer
\eeq

{\it Address}: Max Planck Institute f\"ur Mathematik, Vivatsgasse 7,  Bonn 53111, Germany.
\\
{\it E-mail}: kerner@mpim-bonn.mpg.de
\end{document}